\documentclass[11pt]{article}
\pdfoutput=1
\usepackage{array}
\usepackage{cite}
\usepackage{calc}
\usepackage{caption}
\usepackage{graphicx}
\usepackage{epstopdf}
\usepackage{psfrag}
\usepackage{subcaption}
\usepackage{stfloats}
\usepackage{longtable}
\usepackage{amsmath,amsfonts,amssymb,amsbsy,bm,paralist,theorem,ifthen,color}
\usepackage{algorithm}
\usepackage{algpseudocode}
\usepackage[dvipsnames]{xcolor}
\DeclareMathOperator{\argmin}{\mbox{argmin}}

\def\argmin{\mathop{\rm argmin}}

\def\dom{\mbox{dom\,}}

\def\prox{\mbox{prox}}

\def\QED{\hfill{\bf Q.E.D.}\smskip}

\newcommand{\beq}{\begin{equation}}
\newcommand{\eeq}{\end{equation}}
\newcommand{\st}{{\rm s.t.}}

\bigskip

\newtheorem{theorem}{Theorem}
\newtheorem{lemma}{Lemma}

\newtheorem{remark}{Remark}

\def\smskip{\par\vskip 5 pt}

\def\QED{\hfill{\bf Q.E.D.}\smskip}

\newcommand{\cJ}{\mathcal{J}}

\newcommand{\cF}{\mathcal{F}}

\begin{document}
	\title{Zeroth Order Nonconvex Multi-Agent Optimization over Networks}
		\author{Davood Hajinezhad$^{1}$, Mingyi Hong$^{2}$,  Alfredo Garcia$^{3}$
		\thanks{$^{1}$ SAS Institute, Cary, NC, USA,
			{\tt\small davood.hajinezhad@sas.com}}
		\thanks{$^{2}$  Department of Electrical and Computer Engineering, University of Minnesota, USA, {\tt\small mhong@umn.edu}}
		\thanks{$^{3}$ Department of Industrial and System Engineering, Texas A\&M University, USA,
			{\tt\small alfredo.garcia@exchange.tamu.edu}}}
			\maketitle
	
	\begin{abstract}
		In this paper we consider distributed optimization problems over a multi-agent network, where each agent can only partially evaluate the objective function, and it is allowed to exchange messages with its immediate neighbors. Differently from all existing works on distributed optimization,  our focus is given to optimizing a class of {\it non-convex} problems, and under the challenging setting where each agent can only access the {\it zeroth-order information} (i.e., the functional values) of its local functions. For different types of network topologies such as undirected connected networks or star networks, we develop efficient distributed algorithms and rigorously analyze their convergence and rate of convergence (to the set of stationary solutions). Numerical results are provided to demonstrate the efficiency of the proposed algorithms. 
	\end{abstract}
	
\section{INTRODUCTION}\label{sub:intro} 
Distributed optimization and control has found wide range of applications in emerging research areas such as data-intensive optimization \cite{hong15busmm_spm}, signal and information processing \cite{giannakis2015decentralized}, multi-agent network resource allocation \cite{Tychogiorgosadhoc2013}, communication networks \cite {liao2017distributed}, just to name a few. Typically this type of problems is expressed as minimizing the sum of additively cost functions, given below
\begin{align}\label{eq:sum}
&\min_{x \in \mathbb{R}^M} \; g(x):=\sum_{i=1}^{N} f_i(x),
\end{align}
where $N$ denotes  the number of agents in the network;  $f_i: \mathbb{R}^M\to \mathbb{R}$ represents some (possibly nonsmooth and nonconvex) cost function related to the agent $i$. It is usually assumed that each agent $i$ has complete information on $f_i$, and they can only communicate with their neighbors. Therefore the key objectives of the individual agents are: 1) to achieve consensus with its neighbors about the optimization variable; 2) to optimize the global objective function $g(x)$.

Extensive research has been done on consensus based distributed optimization, but these works are mostly restricted to the family of {\it convex} problems where $f_i(x)$'s are all convex functions.  In \cite{Nedic09subgradient} a first-order method based on the average consensus termed decentralized subgradient (DSG) has been proposed. Following this work, many other first-order algorithms have been proposed to solve distributed convex optimization problems  under different assumptions on the underlying problem. For example in \cite{Nedic09subgradient} DSG is extended to the case where quantized information is used. In \cite{srivastava2011distributed} a local constraint set is added to each local optimization problem. A dual averaging subgradient method is developed and analyzed in \cite{duchi2012dual}. {  In \cite{nedic2015distributed} an algorithm termed subgradient-push has been developed for a time-varying directed network. Other related algorithms can be found in \cite{lobel2011distributed,shi2015proximal,serhat_nips16}. }
 The methods presented so far only converge to a neighborhood of solution set unless using diminishing stepsizes, however using diminishing stepsizes often makes the convergence slow. 
 In order to overcome such a difficulty, recently the authors of \cite{gurbuzbalaban2017convergence} and \cite{shi2014extra} have proposed two methods, named incremental aggregated gradient (IAG)  and  exact first-order algorithm (EXTRA), both of which are capable of achieving fast convergence using constant stepsizes.
 Another class of algorithms for solving problem \eqref{eq:sum} in the convex cases are designed based on primal-dual methods, such as the  Alternating Direction Method of Multipliers (ADMM) \cite{hong14nonconvex_admm_siam,aybat2018distributed}, many of its variants \cite{hong2017stochastic,mokhtari2016dqm}, and distributed dual decomposition method \cite{Terelius201111245}.  

Despite the fact that distributed optimization in convex setting has a broad applicability, many important applications are inherently nonconvex. For example, the resource allocation  in ad-hoc network \cite{Tychogiorgosadhoc2013},  flow control in communication networks \cite{sun2016distributed}, and distributed matrix factorization \cite{hong2016decomposing}, just to name a few.  Unfortunately, without the key assumption of the convexity of $f_i$'s, the existing algorithms and analysis for convex problems are no longer applicable. Recently a few works have started to consider algorithms for nonconvex distributed optimization problems. 
For example, in \cite{zhu2010approximate}  an algorithm based on dual subgradient method has been proposed, but it relaxes the exact consensus constraint. In \cite{bianchi2013convergence} a distributed stochastic projection algorithm has been proposed, and the algorithm converges to KKT solutions when certain diminishing stepsizes are used.  The authors of \cite{hong14nonconvex_admm_siam} proposed an ADMM based algorithm, and they provided one of the first global convergence rate analysis for distributed nonconvex optimization.  More recently, a new convexification-decomposition based approach named NEXT has been proposed in \cite{Lorenzo16}, which utilizes the technique of {\it gradient tracking} to effectively propagate the information about the local functions over the network. In  \cite{hong2016decomposing,hong_prox_pda} a primal-dual based algorithm with provable convergence rate have been designed for distributed nonconvex optimization problem. In \cite{hajinezhad2015nonconvex,davood_NIPS2016_nestt}  the authors proposed primal-dual algorithms for nonconvex optimization problems over a particular network with a central controller.

A key feature for all the above mentioned algorithms, convex or nonconvex, is that they require at least first-order gradient information, and sometime even the second or higher order  information, in order to guarantee global convergence. Unfortunately, in many real-world problems, obtaining such information can be very expensive, if not impossible. For example, in simulation-based optimization \cite{spall2003simulation},  the objective function of the problem under consideration can only be evaluated using repeated simulation.  In certain scenarios of training deep neural network \cite{liannips2016comprehensive},  the relationship between the decision variables and the objective function is too complicated to derive explicit form of the gradient.  Further, in bandit optimization  \cite{agarwal2010optimal}, a player tries to minimize a sequence of loss functions generated by an adversary, and such loss function can only be observed at those points in which the function is realized. 
In these scenarios, one has to utilize techniques from derivative-free optimization, or optimization using zeroth-order information \cite{spall2003simulation}.   Accurately estimating a gradient often requires extensive simulation. In certain application domains, the complexity of each simulation may require significant computational time (e.g. hours). Even when such simulations are parallelized approaches based upon a centralized gradient estimation are impractical due to the need for synchronization; see \cite{fu2015}. In contrast, a zeroth-order distributed approach requires limited simulations for each node and does not need synchronization. 

Recently, Nesterov \cite{nesterov2011random} has proposed a general framework of zeroth-order gradient based algorithms, for both convex and nonconvex problems. It has been shown that for convex (resp. nonconvex) smooth problems the proposed algorithms require  $\mathcal{O}(\frac{M}{\epsilon^2})$   iterations ($M$ denotes the dimension of the problem) to achieve an $\epsilon$-optimal (resp. $\epsilon$-stationary i.e. $\|\nabla f(x)\|^2\le \epsilon$) solution. Further, for both convex and nonconvex problems, the convergence rate for zeroth-order gradient based-algorithms is   at most $\mathcal{O}(M)$ times worse than that of the first-order gradient-based algorithms.  
Ghadimi and Lan \cite{ghadimi2013stochastic} developed a stochastic zeroth-order gradient method which works for convex and nonconvex optimization problems.   Duchi et al. \cite{duchi2015optimal} proposed a stochastic zeroth-order Mirror Descent based  algorithm for solving  stochastic convex optimization problems. In \cite{Gao18} a zeroth-order ADMM algorithm has been proposed for solving convex optimization problems. The complexity of $\mathcal{O}(\frac{1}{\sqrt{T}})$ has been proved for the proposed algorithm,  where $T$ denotes the total number of iterations. Recently an asynchronous stochastic zeroth-order gradient descent
(ASZD) algorithm is proposed in \cite{liannips2016comprehensive} for solving stochastic nonconvex optimization problem. 
\begin{figure}
	\centering
	\hspace{-1cm}\includegraphics[width=1.2 in]{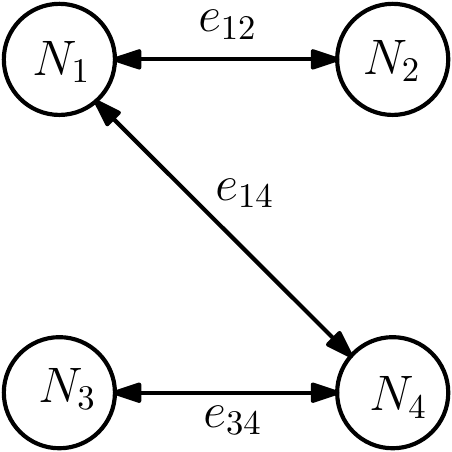}
	\hspace{4cm}\includegraphics[width=1.2 in]{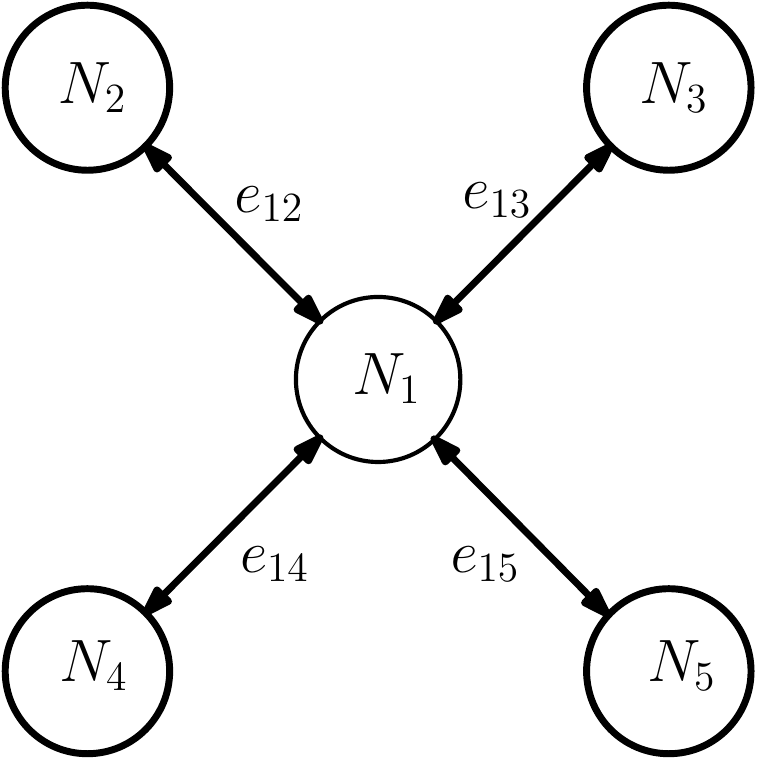}
	\hspace{-.5cm}\caption{ Left: Mesh Network (MNet); \hspace{.2cm}Right: Star Network (SNet)}\label{fig:mnet}
\end{figure}

In this  work we are interested in developing algorithms for the challenging problem of nonconvex distributed optimization, under the setting where each agent $i$ can only access the {\it zeroth-order information} of its local functions $f_i$. For two different types of network topologies, namely, the undirected mesh network (MNet)  (cf. Fig. \ref{fig:mnet}) and the star networks (SNet) (cf. Fig. \ref{fig:mnet}), we develop efficient distributed algorithms and rigorously analyze their convergence and rate of convergence (to the set of stationary solutions). 

In particular, the MNet refers to a network whose nodes are connected to a subset of nodes through an undirected link, and such a network is very popular in applications such as distributed machine learning \cite{forero2011distributed, mateos2010distributed}, and distributed signal processing \cite{Giannakis15, Schizas08}. On the other hand, the SNet has a central controller, which is connected to all the rest of the nodes. Such a network is popular in parallel computing; see for example \cite{icml2015_zhangb15, li2013distributed, hajinezhad2015nonconvex}.  
The main contributions of our work is given below:
\begin{itemize}
	\item For MNet, we design an algorithm capable of dealing with nonconvexity and zeroth-order information in the distributed setting. The proposed algorithm is based upon a primal-dual based zeroth-order scheme, which is shown to converge to the set of stationary solutions of problem \eqref{eq:sum} (with nonconvex but smooth $f_i$'s), in a globally  sublinear manner. \footnote{ Meaning the algorithm converges toward stationary solutions starting from an {\it arbitrary } initial solution.}
	\item For SNet we propose a stochastic primal-dual based method, which is able to further utilize the special structure of the network (i.e., the presence of the central controller) and deal with problem \eqref{eq:sum} with nonsmooth objective. 
	Theoretically, we show that the proposed algorithm also converges to the set of stationary solutions in a globally sublinearly manner.
\end{itemize}
To the best of our knowledge, these algorithms are the first ones for distributed nonconvex optimization that are capable of utilizing zeroth-order information, while possessing global convergence rate guarantees. 

\noindent {\bf Notation.} 
We use $\|a\|$ to denote the Euclidean norm of a vector $a$,  and use $\|A\|$ to denote the spectral norm of matrix $A$. For matrix $A$, $A^\top$ represent its transpose. For a given vector $a$ and matrix $H$, we define $\|a\|^2_H:=a^THa$. The notation {$\langle a,b\rangle$} is used to denote the inner product of two vectors $a$, $b$. To denote an $M\times M$ identity matrix we use $I_M$.  $\mathbb{E}[\cdot]$ denotes taking expectation with respect to all random variables, and $\mathbb{E}_v[\cdot]$ denote taking expectation with respect to the random variable $v$.

\noindent{ \bf Preliminaries.}
We present  some basic concepts and key properties related to derivative-free optimization \cite{nesterov2011random}. 
	Suppose $\mu>0$ is the so-called smoothing parameter, then for a standard Gaussian random vector $\phi\in \mathbb{R}^Q$  the smoothed version of function $\psi:\mathbb{R}^Q\to\mathbb{R}$ is defined as follows 
	\begin{align*}
	\psi_{\mu}(z) = \mathbb{E}_\phi[\psi(z+\mu \phi)]=\frac{1}{(2 \pi)^{\frac{Q}{2}}} \int \psi(z+\mu \phi) e^{-\frac{1}{2}\|\phi\|^2} \,d\phi.
	\end{align*}
	Let us assume that function $\psi$ is $\hat{L}$-smooth (denoted as $\psi\in {\cal C}^{1}_{\hat{L}}$), i.e. there exists a constant $\hat{L}>0$ such that
	\begin{align}\label{eq:def:lsmooth}
	\|\nabla \psi(z_1)-\nabla \psi(z_2)\|\le \hat{L}\|z_1-z_2\|, ~ {  \forall ~ z_1, z_2 \in \mbox{dom}(\psi)}.
	\end{align}
	Then it can be shown that the function $\psi_{\mu} \in \mathcal{C}_{L_\mu}^1$ for some $L_{\mu} \le \hat{L}$, and its gradient is given by Eq. (22) in \cite{nesterov2011random}
	\begin{align} \label{eq:grad_smooth_ver}
	\nabla \psi_{\mu}(z) = \frac{1}{(2 \pi)^{\frac{Q}{2}}} \int \frac{\psi(z+\mu \phi)-\psi(z)}{\mu} \phi e^{-\frac{1}{2}\|\phi\|^2} \,d\phi.
	\end{align}
	Further, for any {$z \in \mathbb{R}^Q$}, it is proved in   \cite[Theorem 1, Lemma 3]{nesterov2011random} that
	\begin{align}
	&{ |\psi_{\mu}(z)-\psi(z)|} \le \frac{\mu^2}{2} \hat{L} Q, \label{eq:diff:g:bd}\\
	&\|\nabla \psi_{\mu}(z) - \nabla \psi(z)\| \le \frac{\mu}{2}\hat{L} (Q+3)^{\frac{3}{2}}, ~\forall z\in\dom(\psi). \label{eq:diff:grad:bd}
	\end{align}
	A stochastic zeroth-order oracle ($\mathcal{SZO}$) takes $z\in \dom(\psi)$ and returns a noisy functional value of $\psi(z)$, denoted by $\mathcal{H}(z; \xi)$, where $\xi\in\mathbb{R}$ is a random variable characterizing the stochasticity of $\mathcal{H}$.  We make the following assumption regarding $\mathcal{H}(z; \xi)$ and $\nabla \psi(z)$.
	
	\noindent{\bf Assumption A.} We assume the following
	\begin{itemize}
		\item [A1.]   $\text{Dom} (\psi)$ is an open set, and  there exists  $K\ge 0$ such that $\forall~ z\in \dom(\psi)$, we have: $\|\nabla \psi(z)\|\le K$;
		\item [A2.] For all $z\in \dom(\psi)$, $\mathbb{E}_\xi[\mathcal{H}(z,\xi)] = \psi(z)$;
		\item [A3.] For some $\sigma\ge 0$, $\mathbb{E}[\|\nabla\mathcal{H}(z,\xi)-\nabla \psi(z)\|^2]\le \sigma^2$, where $\nabla\mathcal{H}(z,\xi)$ denotes any stochastic estimator for $\nabla \psi(z)$. 
	\end{itemize}
	These assumptions are standard in zeroth-order optimization. See for example \cite[Def 1.3 and Lemma 4.2]{Gao18}, \cite[Eq. 4]{nesterov2011random}, and \cite[A3]{ghadimi2013stochastic}.
	Utilizing the $\mathcal{SZO}$ to obtain the functional values, one can show that the following quantity is an unbiased estimator for $\nabla \psi_{\mu}(z)$
		\begin{align}\label{eq:grad:estimate}
		G_\mu(z, \phi, \xi) = \frac{\mathcal{H}(z+\mu\phi, \xi) - \mathcal{H}(z, \xi)}{\mu}\phi,
		\end{align}
		where the constant  $\mu>0$ is a smoothing parameter;  $\phi\in \mathbb{R}^Q$ is a standard  Gaussian random vector.  In particular, we have 
		\begin{align} \label{eq:def_ex_G_mu}
		\mathbb{E}_{\xi, \phi}[G_{\mu}(z, \phi, \xi)] & = \mathbb{E}_\phi \bigg[ \mathbb{E}_{\xi}[G_{\mu}(z, \phi, \xi)\mid \phi] \bigg] = \nabla \psi_{\mu}(z).
		\end{align}
		Furthermore, for given $J$ independent samples of $\{(\phi_j,\xi_j)\}_{j=1}^J$, we define $\bar{G}_{\mu}(z,\phi,\xi)$ as the sample average:  
		\begin{align}\label{eq:J}
		\bar{G}_{\mu}(z,\phi, \xi):=\frac{1}{J}\sum_{j=1}^J G_\mu(z, \phi_j, \xi_j),
		\end{align}
		where $\phi:=\{\phi_j\}_{j=1}^J$, $\xi:=\{\xi_j\}_{j=1}^J$. 
		It is easy to see that for any $J\ge 1$, { $\bar{G}_{\mu}(z,\phi, \xi)$} is also an unbiased estimator of $\nabla \psi_\mu (z)$. Utilizing the above notations and definitions we have the following lemma regarding the  $\bar{G}_{\mu}(z,\phi,\xi)$.
		\begin{lemma}\label{lem:var_bd}\cite[Lemma 4.2]{Gao18}
			Suppose that  Assumption A holds true for function $\psi:\mathbb{R}^Q\to\mathbb{R}$.  Then we have the following
			\begin{align}\label{eq:var_bd}
			\mathbb{E}_{\xi,\phi}[\|{ \bar{G}_{\mu}(z,\phi, \xi)}-\nabla \psi_{\mu}(z)\|^2]\le \frac{\tilde{\sigma}^2}{J},
			\end{align}
			where $\tilde{\sigma}^2:=2Q[K^2+\sigma^2+\mu^2\hat{L}^2Q]$.
		\end{lemma}
 
\section{Zeroth-Order Algorithm over MNet} \label{sec:zos}
\subsection{System Model}
Consider a network of agents represented by a graph $\mathcal{G}:= \{\mathcal{V}, \mathcal{E}\}$, with $|\mathcal{V}| = N$ ($N$ nodes) and $|\mathcal{E}| =E$ ($E$ edges). Each node $v\in\mathcal{V}$ represents an agent in the network, and each edge $e_{ij}=(i,j)\in \mathcal{E}$ indicates that node $i$ and $j$ are neighbors.  Let $\mathcal{N}_i :=\{j\mid (i,j)\in{  \mathcal{E}}\}$ denote the set of neighbors of agent $i$, and assume that $|\mathcal{N}_i|=d_i$.  We assumed that each node can only communicate with its $d_i$ single-hop neighbors in $\mathcal{N}_i$.   

We consider the following reformulation of problem \eqref{eq:sum} 
\begin{align}\label{pr-due}
\min_{z_i\in \mathbb{R}^M}~ \sum_{i=1}^{N}f_i(z_i),\quad\st ~ z_i = z_j, \; \forall~e_{ij}\in\mathcal{E},
\end{align}
where for each agent $i=1,\cdots N$ we introduce a local variable $z_i\in\mathbb{R}^{M}$. If the graph $\mathcal{G}$ is a connected graph, then problem \eqref{pr-due} is equivalent to problem \eqref{eq:sum}.  
For simplicity of presentation let us set ${  Q:=NM}$, and  define a new variable $z:=\{z_i\}_{i=1}^N\in \mathbb{R}^{Q\times1}$.  Throughout this section, we will assume that each function $f_i:\mathbb{R}^{M}\to \mathbb{R}$ is a nonconvex and smooth function. 
Below we present a few network related quantities to be used shortly.
\begin{itemize}
	\item {\it The incidence matrix:} For a given graph $\mathcal{G}$, the {\it incidence matrix} $\tilde{A}\in \mathbb{R}^{E\times N}$ is a matrix where  for each edge $k=(i,j)\in \mathcal{E}$ and when $j> i$, we set $\tilde{A}(k,i)=1$ and $\tilde{A}(k,j)=-1$.  The rest of the entries of $\tilde{A}$ are all zero. For example, for the network in Fig. \ref{fig:mnet} the edge set is ${  \mathcal{E}}=\{e_{12}, e_{14}, e_{34}\}$, therefore the incidence matrix is given by	
	\begin{align*}
	\tilde{A}=\begin{bmatrix}
	1 & -1 &0  &0 \\ 
	1&0  &0  &-1 \\  
	0& 0 &1 &-1
	\end{bmatrix}.	
	\end{align*}
	
	Define the {\it extended incidence matrix} as 
	\begin{align}\label{eq:extended:incidence}
	A:=\tilde{A}\otimes I_M\in\mathbb{R}^{EM\times Q}.
	\end{align}
	\item {\it The degree matrix}:  For a given graph $\mathcal{G}$, the {\it degree matrix} $\tilde{D}\in\mathbb{R}^{N\times N}$ is a diagonal matrix with $\tilde{D}(i,i)= d_i$ where $d_i$ is the degree of node $i$; let {$D :=  \tilde{D}\otimes I_M\in\mathbb{R}^{Q\times Q}$}.	
	\item  {\it The signed/signless Laplacian matrix}: For a given graph $\mathcal{G}$ with its extended incidence matrix given by \eqref{eq:extended:incidence}, its signed and signless Laplacian matrices are expressed as
	\begin{subequations}
	\begin{align}
	L^{-} &:=A^\top A\in\mathbb{R}^{Q\times Q},\label{eq:slap}\\
	L^{+}&:=2D-A^\top A\in\mathbb{R}^{Q\times Q},\label{eq:sllap}
	\end{align}	
	\end{subequations}
\end{itemize}  
respectively. Using the above notations, one can easily check that problem \eqref{pr-due} can be written compactly as below 
\begin{align}\label{pro:compact}
\min_{z\in \mathbb{R}^Q} \; g(z)=\sum_{i=1}^{N}f_i(z_i),\quad \st \; \;  Az = 0,
\end{align}
where we have defined $z:=\{z_i\}_{i=1}^N\in \mathbb{R}^{Q\times1}$.  
The Lagrangian function for this problem is defined by
\begin{equation}\label{eq:lag1}
L(z,\lambda) := g(z)+\langle\lambda,Az\rangle,
\end{equation}
where {$\lambda \in\mathbb{R}^{EM\times 1}$} is the dual variable  associated with the constraint $Az = 0$.  The stationary solution set for the problem \eqref{pro:compact} is given by
\begin{align}\label{eq:stationary}
S = \{(z^*, \lambda^*) \mid \nabla_z L(z^*,\lambda^*)=0 ~\text{and}~ Az^*=0\},
\end{align}
where $\nabla_z L(z^*,\lambda^*)$ denotes the gradient of Lagrangian function with respect to the variable $z$ evaluated at $(z^*,\lambda^*)$.
\subsection{The Proposed Algorithm}
In this subsection we present a {\underline{Z}eroth-\underline{O}rder {\underline N}onconv\underline{E}x, over \underline{M}Net (ZONE-M)} algorithm which is capable of solving distributed nonconvex optimization problem in an efficient manner [towards approximating the  stationary solution as defined in \eqref{eq:stationary}].  To proceed, let us first construct the augmented Lagrangian (AL) function for problem \eqref{pro:compact} 
\begin{equation}\label{eq:lag}
L_\rho(z,\lambda) := g(z)+\langle\lambda,Az\rangle+ \frac{\rho}{2}\big\|Az\big\|^2,
\end{equation}
where {$\lambda \in\mathbb{R}^{EM\times 1}$} is the dual variable  associated with the constraint $Az = 0$, and $\rho>0$ denotes the  penalty parameter. To update the primal variable $z$, the AL is first approximated using a quadratic function with a degree-matrix weighted proximal term $\|z-z^r\|^2_D$, followed by one step of zeroth-order gradient update to optimize such a quadratic approximation. After the primal update, an approximated dual ascent step is performed to update $\lambda$. The algorithm steps are detailed in Algorithm \ref{alg:zenith}.  { In this algorithm, $\mathcal{H}_i(z,\xi)$ denotes a noisy functional value returned by a $\mathcal{SZO}$ associated to the local function $f_i$, and we assumed that $\mathcal{H}_i(z,\xi)$ satisfies Assumption A for all $i=1,\cdots, N$.} Note that the ZONE-M  is a variant of the popular method called {\it Method of Multipliers} (MM), whose steps are expressed below \cite{Hestenes} 
\begin{align}
z^{r+1}&=\argmin_{z\in\mathbb{R}^Q} L_\rho(z, \lambda^r),\label{eq:aug:z}\\
\lambda^{r+1} &= \lambda^{r}+\rho Az^{r+1}.
\end{align}
However, for the problem that is of interest in this paper, the MM method is not applicable because of the following reasons: 1) The optimization problem \eqref{eq:aug:z} is not easily solvable to global optima because it is nonconvex, and we only have access to zeroth-order information; 2) It is not clear how to implement the algorithm in a distributed manner over the MNet.  
In contrast, the primal step of the ZONE-M algorithm \eqref{eq:z:zero} utilizes zeroth-order information and can be performed in closed-form. Further, as we  elaborate below, combining the primal and the dual steps of ZONE-M yields a fully distributed algorithm. 

\begin{algorithm}[tb]
	\caption{The ZONE-M Algorithm}
	\label{alg:zenith}
	\begin{algorithmic}[1]
		\State {\bfseries Input:} ${z}^0\in\mathbb{R}^{Q}$, {$\lambda^0=0_{EM}$}, $D\in\mathbb{R}^{Q\times Q}$, $A\in\mathbb{R}^{EM\times Q}$, $T\ge 1$, $ J\ge1$, { $\mu>0$}
		\For{$r=0$ {\bfseries to} ${  T-1}$}  
		\Statex For each $i=1,\cdots, N$, 
		 generate $\phi^r_{i,j}\in\mathbb{R}^M$, $j=1,2,\cdots,J$  from an i.i.d standard Gaussian distribution and
		calculate $\bar{G}_{\mu,i}(z_i^r,\phi^r_{i},\xi^r_{i})\in\mathbb{R}^{M}$ by
		\begin{align}\label{eq:grad:approx}
		&{ \bar{G}_{\mu,i}(z_i^r,\phi^r_{i},\xi^r_{i})}=\frac{1}{J}\sum_{j=1}^J\frac{\mathcal{H}_i(z_i^r+\mu\phi^r_{i,j}, \xi^r_{i,j}) - \mathcal{H}_i(z_i^r, \xi^r_{i,j})}{\mu}\phi^r_{i,j},
	\end{align}
		where we have defined  $\phi_i^r:=\{\phi^r_{i,j}\}_{j=1}^{J}$,  $\xi_i^r:=\{\xi^r_{i,j}\}_{j=1}^{J}$; Define ${  G^{J,r}_{\mu}}:=\{\bar{G}_{\mu,i}(z_i^r,\phi^r_{i},\xi^r_{i})\}_{i=1}^N\in\mathbb{R}^{Q}$.
		\Statex Update $z$ and $\lambda$ by 
		\begin{align}
		z^{r+1}& = \argmin_{z} ~\langle {G^{J,r}_{\mu}}+A^\top\lambda^r+\rho A^\top Az^r,z-z^r\rangle +\rho\|z-z^r\|^2_{D},\label{eq:z:zero}\\
		\lambda^{r+1}& = \lambda^r +\rho A z^{r+1}.\label{eq:lambda:zero}
		\end{align}
		\EndFor
		\State { Choose uniformly randomly $u\in\{0,1,\cdots,T-1\}$}
		\State {\bfseries Output:} $(z^u, \lambda^u)$.
	\end{algorithmic}
\end{algorithm}
 
To illustrate the distributed implementation of the proposed method, let us transform the ZONE-M algorithm to a {\it primal only} form. To this end, let us write down the optimality condition for \eqref{eq:z:zero} as 
\begin{align}\label{eq:opt:cond}
{  G^{J,r}_{\mu}}+A^\top\lambda^r+\rho A^\top Az^r+2\rho D(z^{r+1}-z^r)=0.
\end{align}
Utilizing the definitions in \eqref{eq:slap}, and \eqref{eq:sllap}, we have the following identity from \eqref{eq:opt:cond}
\begin{align}\label{eq:kkt:z:compact}
{  G^{J,r}_{\mu}}+A^\top\lambda^r+2\rho D z^{r+1}-\rho L^{+}z^r=0.
\end{align}
Let us replace $r$ in equation \eqref{eq:kkt:z:compact} with $r-1$ to get
\begin{align}\label{eq:kkt:x:compact:r-1}
{  G^{J,r-1}_{\mu}}+A^\top\lambda^{r-1}+2\rho D z^{r}-\rho L^{+}z^{r-1}=0.
\end{align}
Now rearranging the terms in \eqref{eq:lambda:zero} and using the definition in \eqref{eq:slap} we have 
\begin{align}\label{eq:due:arng}
A^\top (\lambda^{r}-\lambda^{r-1})=\rho A^\top A z^r = \rho L^{-}z^r. 
\end{align}
Subtracting equation \eqref{eq:kkt:x:compact:r-1} from \eqref{eq:kkt:z:compact} and utilizing \eqref{eq:due:arng} yield
\begin{align*}
&{  G^{J,r}_{\mu} - G^{J,r-1}_{\mu}} + \rho L^{-} z^r
+ 2\rho D (z^{r+1}-z^r)-\rho L^{+}(z^{r}-z^{r-1}) = 0.
\end{align*}
Rearranging terms in the above identity, we obtain
\begin{align}\label{eq:z:compact}
z^{r+1} &= z^r -\frac{1}{2\rho} D^{-1}\bigg[{  G^{J,r}_{\mu} - G^{J,r-1}_{\mu}}\bigg]
+ \frac{1}{2} D^{-1} (L^{+}-L^{-})z^r-\frac{1}{2} D^{-1} L^{+}z^{r-1}.
\end{align} To implement such iteration, it is easy to check (by utilizing the definition of $L^{+}$ and $L^{-}$) that each agent $i$ performs the following local computation 
\begin{align}\label{eq:zenith:distributed}
z_i^{r+1}& = z_i^r -\frac{1}{2\rho d_i} \bigg[\bar{G}_{\mu,i}(z_i^{r}, \phi_i^{r},\xi_i^{r}) - \bar{G}_{\mu,i}(z_i^{r-1}, \phi_i^{r-1},\xi_i^{r-1})\bigg] \nonumber\\
&+ \sum_{j\in \mathcal{N}_i}\frac{1}{d_i} z_j^r -\frac{1}{2}\bigg(\sum_{j\in \mathcal{N}_i}\frac{1}{d_i} z_j^{r-1} + z_i^{r-1}\bigg),
\end{align} where $\bar{G}_{\mu,i}(z_i^{r}, \phi_i^{r},\xi_i^{r})$ is defined in \eqref{eq:grad:approx}. Clearly, this is a fully decentralized algorithm, because to carry out such an iteration, each agent $i$ only requires the knowledge about its local function [i.e., $\bar{G}_{\mu,i}(z_i^{r}, \phi_i^{r},\xi_i^{r})$,  ${\bar{G}_{\mu,i}(z_i^{r-1}, \phi^{r-1}_i,\xi^{r-1}_i)}$, $z^r_i$ and $z^{r-1}_i$], as well as  information from the agents in its neighborhood $\mathcal{N}_i$. 
	\begin{remark}
	The single variable iteration derived in \eqref{eq:z:compact} takes a similar form as the EXTRA algorithm proposed in \cite{shi2014extra}, which uses the first-order gradient information. In EXTRA, the iteration is given by (for $r\geq 2$)
		\begin{align*}
		z^{r+1} -z^r &= Wz^r  -\frac{I_Q+W}{2}z^{r-1}-\alpha \bigg[\nabla g(z^r) - \nabla g(z^{r-1})\bigg],
		\end{align*} where $W$ is a double stochastic matrix. 
	
	In the ZONE-M algorithm, let us define $W: = \frac{1}{2}D^{-1} (L^{+}-L^{-})$, which is a row stochastic matrix. Then iteration \eqref{eq:z:compact} becomes 
	\begin{align*}
	z^{r+1}-z^r &= Wz^r  -\frac{I_Q+W}{2}z^{r-1}-\frac{1}{2\rho} D^{-1}\bigg[G^{J,r}_{\mu} - G^{J,r-1}_{\mu}\bigg],
	\end{align*} which is similar to the EXTRA algorithm. The key difference is that our algorithm { utilizes} {\it zeroth-order} information, to deal with {\it nonconvex} problems, while the EXTRA algorithm requires {\it first-order} (gradient) information, and it only deals with {\it convex} problems. 
\end{remark}
\subsection{The Convergence Analysis of ZONE-M}\label{sec:analysis}
In this subsection we provide the convergence analysis for the ZONE-M algorithm.  We  make the following  assumptions. 

\noindent{\bf  Assumptions  B.} 
\begin{itemize}
	\item [B1.] For all $i\in\{1,2,\cdots,N\}$,  $\text{Dom} (f_i)$ is an open set, and  there exists  $K_i\ge 0$ such that $\forall~ z_i\in \dom(f_i)$, we have: $\|\nabla f_i(z_i)\|\le K_i$;
	\item [B2.] For all $z_i\in \dom(f_i)$, $\mathbb{E}_{\xi_i}[\mathcal{H}_i(z_i,\xi_i)] = f_i(z_i)$;
	\item [B3.] For some $\sigma_i\ge 0$, $\mathbb{E}[\|\nabla\mathcal{H}_i(z_i,\xi_i)-\nabla f_i(z_i)\|^2]\le \sigma_i^2$, where $\nabla\mathcal{H}_i(z_i,\xi_i)$ denotes any stochastic estimator for $\nabla f_i(z_i)$. 
	\item[B4.]  Function $g:=\sum_{i=1}^{N}f_i$ is $\hat{L}$-smooth, which satisfies \eqref {eq:def:lsmooth}.
	\item[B5.] There exists a constant $\delta>0$ such that
	\begin{align}\label{eq:delta_bd}
	   \exists ~\underline{g}>-\infty, \quad \st\quad g(z)+\frac{\delta}{2}\|Az\|^2\ge \underline{g},\; \forall~z\in\mathbb{R}^Q. 
    \end{align}
\end{itemize}
 Without loss of generality we can set $\underline{g}=0$. 
A few examples of nonconvex functions that satisfy the Assumptions B are provided below:
\begin{itemize}
	\item The sigmoid function $\mbox{sig}(z)=\frac{1}{1+e^{-z}}$
	\item The function $\tanh(z)=\frac{1-e^{-2z}}{1+e^{-2z}}$
	\item The function $2\mbox{logit}(z) =\frac{2e^{z}}{e^z+1}=1+ \tanh(\frac{z}{2})$
\end{itemize} 
 Let us define the gradient of smoothed version of function $g$ denoted by $\nabla g_\mu$ similar to \eqref{eq:grad_smooth_ver}. From Assumption [B4] and the preliminary results we conclude that $\nabla g_\mu$ is $L_\mu$-smooth, where $L_\mu\leq \hat{{L}}$. Also, one can simply check that whenever all $f_i$'s satisfy Assumptions [B1-B3], the function $g:=\sum_{i=1}^N f_i$ also satieties a similar sets of assumptions as Assumptions [B1-B3]. In particular, there exist constants $K_g$ and $\sigma_g$ such that $\|\nabla g(z)\|\leq K_g$, and $\mathbb{E}[\|\nabla\mathcal{H}(z,\xi)-\nabla g(z)\|^2]\le \sigma_g^2$. As a result, we can apply Lemma \ref{lem:var_bd} for function $g:\mathbb{R}^Q\to\mathbb{R}$. Therefore, setting 
	\begin{align}\label{eq:sigma_g}
	\tilde{\sigma}^2_g:=2Q[K_g^2+\sigma_g^2+\mu^2+\hat{L}^2Q],
	\end{align} 
	we have 
\begin{align}\label{eq:var_bd:mnet}
\mathbb{E}_{\xi,\phi}[\|{{G}^{J,r}_{\mu}}-\nabla g_{\mu}(z)\|^2]\le \frac{\tilde{\sigma}_g^2}{J}.
\end{align} 

Let ${\mathcal{F}^{r}:=\{(\xi^t, \phi^t)\}_{t=1}^{r}}$ be the $\sigma$-field generated by the entire history of algorithm up to iteration $r$. Let $\sigma_{\min}$ be the smallest nonzero eigenvalue of matrix $A^\top A$. Additionally, we define 
$w^r:=(z^{r+1}-z^r) -(z^{r}-z^{r-1}).$  Further, to facilitate the analysis let us list a few  relationships below. 
\begin{itemize}
	\item For any given vectors $a$ and $b$ we have
	\begin{align}
	&\langle b-a, b\rangle = \frac{1}{2}(\|b\|^2+\|a-b\|^2-\|a\|^2), \label{eq:rel1}\\
	& \langle a,b\rangle \le \frac{1}{2\epsilon}\|a\|^2 + \frac{\epsilon}{2}\|b\|^2; \quad \forall~ \epsilon>0. \label{eq:rel3}
	\end{align}
	\item For $n$ given vectors  $a_i$ we have the following 
	\begin{align}\label{eq:rel2}
	\bigg\|\sum_{i=1}^na_i\bigg\|^2\le n\sum_{i=1}^{n}\big\|a_i\big\|^2.
	\end{align}
\end{itemize}

Our convergence analysis consists of the following main steps: First we show that the successive difference of the dual variables, which represents the constraint violation, is bounded by a quantity related to the primal variable. Second we construct a special potential function whose behavior is tractable under a specific parameter selection. Third, we combine the previous results to obtain the main convergence rate analysis. Below we provide a sequence of lemmas and the main theorem. The proofs  are provided in Appendix A. 
Unless otherwise stated, throughout this section the expectations are taken with respect to {$(\xi^{r+1}, \phi^{r+1})$ conditioning on the filtration $\mathcal{F}^{r}$ defined previously.} 

Our first lemma bounds the change of the dual variables (in expectation) by that of the primal variables. This lemma will be used later to control the progress of the dual step of the algorithm. 
	\begin{lemma}\label{lemma:mu:bound}
	  Suppose Assumptions B hold true. Then for   $r\geq 1$ we have the following relation	
	\begin{align}\label{eq:mu:difference:bound}
	\mathbb{E}\|\lambda^{r+1}-\lambda^{r}\|^2&\le {\frac{9\tilde{\sigma}_g^2}{J\sigma_{\min}}}+\frac{6L^2_\mu}{{\sigma_{\min}}}\mathbb{E}	\|z^r-z^{r-1}\|^2+\frac{3\rho^2\|L^{+}\|}{\sigma_{\min}}\mathbb{E}\|w^r\|_{L^{+}}^2,
	\end{align}
\end{lemma} 
where $\tilde{\sigma}_g$ is defined in \eqref{eq:sigma_g}.

To proceed, we need to construct a  potential function so that the behavior of the algorithm can be made tractable. For notational simplicity let us define $L^{r+1}_{\rho}:=L_\rho(z^{r+1},\lambda^{r+1})$.  Also let $c>0$ to be some positive constant (to be specified shortly), and set  {  $k:=2\bigg(\frac{6\hat{L}^2}{\rho\sigma_{\min}}+\frac{3c\hat{L}}{2}\bigg)$}. { Define $V^{r+1}:=\frac{\rho}{2}\bigg(\|Az^{r+1}\|^2+\|z^{r+1}-z^r\|^2_{B}\bigg)$, where $B:=L^++\frac{k}{c\rho}I_Q$.} Using these notations, we define a potential function in the following form
\begin{align}\label{eq:pot_def}
P^{r+1}:=L^{r+1}_{\rho}+ cV^{r+1}.
\end{align}

The following lemma analyzes the behavior of the potential function as the ZONE-M algorithm proceeds.
\begin{lemma}\label{lem:bd:pot}
Suppose Assumptions B hold true, and parameters $c$ and $\rho$ satisfy the following conditions 	
	\begin{align}
	&c>\frac{6\|L^{+}\|}{\sigma_{\min}},\quad\rho>\max\big(\frac{-b+\sqrt{b^2-8d}}{4},\delta,\hat{L}/2\big),\label{eq:bd:c}
	\end{align}
where 
	\begin{align*}
	{  b=-\hat{L}(\hat{L}+4c+1)-3},~ d=-\frac{12\hat{L}^2}{\sigma_{\min}}.
	\end{align*}
	Then for some constants $c_1, c_2, c_3>0$, the following inequality holds true for   $r\geq 1$
	\begin{align}
	\mathbb{E}\bigg[P^{r+1}-P^{r}\bigg]&\le \frac{k-c_1}{2}\mathbb{E}\|z^{r+1} - z^r\|^2 -c_2\mathbb{E}\|w^r\|_{L^{+}}^2 
	+  c_3\frac{\tilde{\sigma}_g^2}{J}+\frac{3\mu^2({Q}+3)^3}{8},\label{eq:diff:bd}
	\end{align}
	where we have defined the following constants
	\begin{align}\label{eq:c}
	&c_1:= 2\rho-\hat{L}^2-(c+1)\hat{L}-3>0,\\
	& c_2:=\bigg(\frac{c \rho}{2}-\frac{3\rho\|L^{+}\|}{\sigma_{\min}}\bigg)>0, ~ c_3:=\frac{9}{\rho\sigma_{\min}}+\frac{3+6c\hat{L}}{2\hat{L^2}}>0\nonumber.
	\end{align}
\end{lemma}
	We can readily observe that using the choice of $c$ in \eqref{eq:bd:c}, $c_2$ is positive. Further for any fixed $c$, it is possible to make $\rho$ sufficiently large such that $k-c_1<0$. Therefore in expectation, the potential function {\it decreases} in $\mathbb{E}[\|z^{r+1}-z^r\|^2]$ and $\mathbb{E}\|w^t\|^2_{L^+}$, while it increases in  constants proportional to $\mu^2$ and $\frac{1}{J}$. Later we will leverage this result by properly choosing $\mu$, and $J$  to derive the convergence rate of the algorithm. 

The key insight obtained from this step is that, a conic combination of augmented Lagrangian function, as well as the constraint violation can {  serve} as the potential function that guides the progress of the algorithm. We expect that such construction is of independent interest. It will be instrumental in analyzing other (probably more general) nonconvex primal-dual type algorithms. 

The next lemma shows  that $P^{r+1}$ is lower bounded.
\begin{lemma}\label{lemma:lower:bound}
	 Suppose Assumptions B hold true, and the constant $c$ is picked large enough such that   
	 \begin{align}
 	c\geq \frac{2\|L^+\|}{\sigma_{\min}}.
	 \end{align}
	Then the statement below holds true
		{
		\begin{align}\label{eq:lower:bound}
		\exists~\underline{P}  \quad \st \quad   \mathbb{E}[P^{r+1}]\ge \underline{P}>-\infty, \quad \forall~r\geq 1.
		\end{align}}
	where $\underline{P} $ is a constant that is independent of total number of iterations $T$. 
\end{lemma}

To present our main convergence theorem,  we need to  measure the gap between the current iterate to the set of stationary solutions.  To this end, consider the following gap function 
\begin{align}\label{eq:opt_gap}
\Phi(z^r,\lambda^{r-1}): = \mathbb{E}\bigg[\|\nabla_{z}L_{\rho}(z^r, \lambda^{r-1})\|^2+\|A z^r\|^2\bigg].
\end{align}
 It can be easily checked that $\|\nabla_{z}L_{\rho}(z^*, \lambda^{*})\|^2+\|A z^*\|^2= 0$ if and only if $(z^*, \lambda^*)$ is a stationary solution of the problem \eqref{pro:compact}. For notational simplicity let us write $\Phi^{r}:=\Phi(z^r, {  \lambda^{r-1}})$.
The result below quantifies the convergence rate of ZONE-M.
\begin{theorem}\label{thm:conv}
	Consider the ZONE-M algorithm. 	Suppose Assumptions B hold true,   the penalty parameter $\rho$ satisfies the condition given in Lemma \ref{lem:bd:pot}, and the constant $c$ satisfies  {$c\geq \frac{6\|L^{+}\|}{\sigma_{\min}}.$}
	Then there exists constants $\gamma_1, \gamma_2, \gamma_3>0$ such that we have the following bound
	\begin{align}
	\mathbb{E}_u[\Phi^u]\le \frac{\gamma_1}{T} + \frac{\gamma_2\tilde{\sigma}_g^2}{J}+\gamma_3\mu^2. 
	\end{align}
 The explicit value for constants $\gamma_1, \gamma_2$, and $\gamma_3$ can be expressed as the following: Let 
 \begin{align*}
 \alpha_1=8\hat{L}+2\rho^2\|L^+\|^2, \; \alpha_2=\frac{6\hat{L}}{\rho^2\sigma_{\min}}, \; \alpha_3=\frac{3\|L^+\|}{\sigma_{\min}},
 \end{align*}
 and $c_1,c_2$ and $c_3$ are constants given in equation \eqref{eq:c}. Let us set $\zeta=\frac{\max(\alpha_1+\alpha_2, \alpha_3)}{\min(\frac{c_1-k}{2}, c_2)}$, { where $c_1:= 2\rho-\hat{L}^2-(c+1)\hat{L}-3$, and $k:=2\bigg(\frac{6\hat{L}^2}{\rho\sigma_{\min}}+\frac{3c\hat{L}}{2}\bigg)$.} Then we have the following expression 
 \begin{align*}
 \gamma_1&=\zeta\mathbb{E}[P^1-\underline{P}]+\alpha_2\mathbb{E}\|z^1-z^0\|^2\\ \gamma_2&=\zeta c_3+\frac{9+4\rho^2\sigma_{\min}}{\rho^2\sigma_{\min}}, ~\gamma_3=\frac{3}{8}\zeta+2\hat{L}^2.
 \end{align*}
\end{theorem}
\begin{remark}
From Theorem \ref{thm:conv} we can observe that the complexity bound of the ZONE-M depends on  $\tilde{\sigma}_g$, and the smoothing parameter $\mu$ . Therefore, no matter how many iterations we run the algorithm, it always converges to a neighborhood of a KKT point, which is expected when only zeroth-order information is available; see  \cite[Theorem 4.4]{Gao18}, and \cite[Theorem 3.2] {ghadimi2013stochastic}. 
 Nevertheless, if we choose $J\in\mathcal{O}(T)$, and $\mu\in\mathcal{O}(\frac{1}{\sqrt{T}})$, we can achieve the following bound
\begin{align}
\mathbb{E}_u[\Phi^u]\le \frac{\gamma_1}{T} +\frac{\gamma_2\tilde{\sigma}_g^2}{T}+\frac{\gamma_3}{T}.
\end{align}
This indicates that ZONE-M converges in a sublinear rate.
\end{remark}
\begin{remark}\label{rem:increas}
 Our  bound   on $\rho$ derived in \eqref{eq:bd:c} can be loose because it is obtained based on the the worst case analysis. In practice one may start with a small $\rho$ and gradually increase it until reaching the theoretical bound. In the numerical experiments, we will see that such a strategy often leads to faster empirical convergence.
\end{remark}

\section{Zeroth-Order Algorithm over SNet}\label{sec:snet}
In this section we focus on multi-agent optimization problem over SNet (cf. Fig. \ref{fig:mnet}). 
 We propose the \underline{Z}eroth-\underline{O}rder {\underline N}onconv\underline{E}x, over \underline{S}Net (ZONE-S) algorithm for the multi-agent optimization problem. 
\subsection{System Model}
Let us consider the following problem

 \begin{align}\label{eq:sum:nonsmoth}
 &\min_{x \in X} \; g(x):=\sum_{i=1}^{N} f_i(x) + r(x),
 \end{align} where $X\subseteq \mathbb{R}^M$ is a closed and convex set, $f_i: \mathbb{R}^M\to \mathbb{R}$ is smooth possibly nonconvex function, and $r: \mathbb{R}^M\to \mathbb{R}$ is a convex possibly nonsmooth function, which is usually used to impose some regularity to the solution. Let us set $f(x): =\sum_{i=1}^{N} f_i(x)$ for notational simplicity. Note that this problem is slightly more general than the one solved in the previous section [i.e., problem \eqref{eq:sum} with smooth objective function], because here we have included constraint set $X$ and the nonsmooth function $r(x)$ as well.

 We note that many first-order algorithms have been developed for solving problem \eqref{eq:sum:nonsmoth}, including SGD \cite{robbins51}, SAGA \cite{Defazio14},  SVRG \cite{Johnson13}, and NESTT \cite{davood_NIPS2016_nestt}, but it is not clear how to {  adapt} these methods and their analysis to the case with non-convex objective and zeroth-order information. 
 
Similar to the  problem over MNet, here we split the variable $x\in \mathbb{R}^M$ into $z_i\in \mathbb{R}^M$, and reformulate problem \eqref{eq:sum:nonsmoth} as
\begin{align}\label{eq:sum:nonsmoth:dist}
\min_{x, z}~  \sum_{i=1}^{N} f_i(z_i) + h(x) ~\st ~ x=z_i, ~ \forall ~i=1,\cdots, N, 
\end{align} where $h(x):=r(x)+\iota_X(x)$, [$\iota_X(x)=0$ if $x\in X$, otherwise $\iota_X(x)=\infty$ ]. 
In this formulation we have assumed that  for $i=1,2,\cdots N$,  $f_i$ is the local function for agent $i$, and $h(x)$ is handled by the central controller. Further, agent $i$ has access to the stochastic functional values of $f_i$ through the $\mathcal{SZO}$ as described in preliminaries. 
\subsection{Proposed Algorithm}
The proposed algorithm is again a primal-dual based scheme. The augmented Lagrangian function for problem \eqref{eq:sum:nonsmoth:dist} is given by
\begin{align*}
&L_\rho\big(z, x; \lambda\big) \nonumber=\sum_{i=1}^{N}\bigg(f_i(z_i)+ \langle \lambda_i, z_i-x\rangle+\frac{\rho_i}{2}\|z_i-x\|^2\bigg)  + h(x),
\end{align*} where $\lambda_i$, and $\rho_i$ are  respectively the dual variable and the penalty parameter associated with the constraint $z_i=x$. Let $\lambda:=\{\lambda_i\}_{i=1}^{N}$, {  $\rho:=\{\rho_i\}_{i=1}^{N}\in \mathbb{R}_{++}^N$}.
To proceed, let us introduce the following function for agent $i$
	\begin{align}\label{eq:al:approx}
	U_{\mu,i}(z_i,x; \lambda_i) &= f_i(x) + \langle \bar{G}_{\mu,i}(x, \phi,\xi), z_i-x\rangle+\langle \lambda_{i}, z_i-x\rangle
	+\frac{\alpha_i\rho_i}{2}\|z_i-x\|^2.
	\end{align}
In the above expression $\alpha_i>0$  is a positive constant, and  $\bar{G}_{\mu,i}(x,\phi,\xi)$ is given by
\begin{align}\label{eq:grad:estimate:zone-ns}
\bar{G}_{\mu,i}(x,\phi,\xi) = \frac{1}{J}\sum_{j=1}^J\frac{\mathcal{H}_i(x+\mu\phi_{j}, \xi_{j}) - \mathcal{H}_i(x, \xi_{j})}{\mu}\phi_{j},
\end{align} where $\mathcal{H}_i(x, \xi)$ is a noisy version of $f_i(x)$ obtained from $\mathcal{SZO}$ { and satisfies Assumption A}, $\mu>0$ is  smoothing parameter, $\phi_{j}\in\mathbb{R}^{M}$ is a standard Gaussian random vector, $\xi_{j}$ represents the noise related to the $\mathcal{SZO}$ output, and we set $\phi=\{\phi_{j}\}_{j=1}^J$, and $\xi=\{\xi_{j}\}_{j=1}^J$. 
  To see more details about the characteristics  of function $U_{\mu,i}(z_i,x; \lambda_i)$ the readers are refereed to \cite{davood_NIPS2016_nestt}.

The proposed algorithm is described below. At the beginning of iteration  $r+1$ the central controller broadcasts $x^r$  to everyone. An agent indexed by $i_r\in\{1,2,\cdots N\}$ is then randomly picked with some probability of $p_{i_r}$, and this agent optimizes $U_{\mu,i_r}(z_i,x^r, \lambda^r)$  [defined in \eqref{eq:al:approx}], and updates its dual variable $\lambda_{i_{r}}$.
The rest of the nodes $j\ne i_r$  simply set $z_j^{r+1}=x^r$, and $\lambda_j^{r+1}=\lambda_j^r$. Finally the central controller updates the variable $x$ by minimizing the augmented Lagrangian. The pseudo-code of the ZONE-S algorithm is presented in Algorithm \ref{alg:zenith-ns}.   
\begin{algorithm}[tb]
	\caption{The ZONE-S Algorithm}
	\label{alg:zenith-ns}
	\begin{algorithmic}[1]
		\State {\bf Input:} ${x}^0 \in \mathbb{R}^M$, $\lambda^0\in \mathbb{R}^M$, $T\ge 1$, $J\ge 1$, $\mu>0$
		\For {$ r=1$ {\bf to} $T$},
		\Statex {\it In central controller}: Pick $i_r$ from $\{1,2,\cdots, N\}$ with probability  $p_{i_r}=\frac{\sqrt{L_{\mu,i_r}}}{\sum_{i=1}^N \sqrt{L_{\mu,i}}}$. Generate	 $\phi_j^r\in\mathbb{R}^M$, $j=1,2,\cdots,J$ from an i.i.d standard Gaussian distribution
		\Statex {\it In agent $i_r$}: Calculate ${ \bar{G}}_{\mu,{i_r}}(x^r,\phi^r,\xi^r)$ using
		\begin{align}
		&\bar{G}_{\mu,{i_r}}(x^r,\phi^r,\xi^r) 
			= \frac{1}{J}\sum_{j=1}^J\frac{\mathcal{H}_{i_r}(x^r+\mu\phi_j^r, \xi_j^r) - \mathcal{H}_{i_r}(x^r, \xi_j^r)}{\mu}\phi_j^r,
		\end{align} 
		where we set $\phi^r=\{\phi^r_j\}_{j=1}^J$, and $\xi^r=\{\xi^r_j\}_{j=1}^J$.
		\Statex {\it In all agents}: Update $z$, and $\lambda$ by
		\begin{align}
		z_{i_r}^{r+1}&=x^r-\frac{1}{\alpha_{i_r}\rho_{i_r}}\bigg[\lambda^r_{i_r}+\bar{G}_{\mu,{i_r}}(x^r,\phi^r,\xi^r)\bigg];\label{eq:z_i:zone-ns}&&\\
		\lambda^{r+1}_{i_r}&=\lambda_{i_r}^{r}+\alpha_{i_r}\rho_{i_r}\bigg(z^{r+1}_{i_r}-x^{r}\bigg); \label{eq:lam:zone-ns}&&\\
		\lambda^{r+1}_{j}&=\lambda_{j}^{r}, \quad z^{r+1}_{j} = x^r, \quad \forall~j\ne {i_r}. \label{eq:z_i:zone-ns2}&&
		\end{align}
		{\it In central controller}: Update $x$ by
		\begin{align}
		\quad \quad  x^{r+1}&=\arg\min_{x\in X}L_{  \rho}(z^{r+1}, x;  \lambda^{r}).\label{eq:x:zone-ns}&&
		\end{align}
		\EndFor
		\State  Choose uniformly randomly $u\in\{1,2,\cdots,T\}$.
		\State{\bf Output:} $x^u$.
	\end{algorithmic}
\end{algorithm} 
\subsection{Convergence Analysis of ZONE-S}
We make the following assumptions in this part. 

 \noindent{\bf Assumption C.}
\begin{itemize}
	\item [C1.] {  $\text{Dom} (f_i)$ is an open set, and  there exists  $K_i\ge 0$ such that $\forall~ x\in \dom(f_i)$, we have: $\|\nabla f_i(x)\|\le K_i$;}
	\item [C2.] For all $x\in \dom(f_i)$, $\mathbb{E}_\xi[\mathcal{H}_i(x,\xi)] = f_i(x)$;
	\item [C3.] For some $\sigma_i\ge 0$, $\mathbb{E}[\|\nabla\mathcal{H}_i(x,\xi)-\nabla f_i(x)\|^2]\le \sigma_i^2$.
	\item [C4.] For $i=1,2,\cdots, N$, function $f_i$ and $f$ are $L_i$-smooth, and $L$-smooth respectively.
	\item [C5.] The function $g(x)$ is bounded from below over $X\cap \hbox{int}(\hbox{dom }(g))$.
	\item [C6.] The function $r(x)$ is convex but possibly nonsmooth. 
\end{itemize}

Let us define $\tilde{\sigma}_i:=2M[K_i^2+\sigma_i^2+\mu^2{L}_i^2M]$, and set $\tilde{\sigma}^2:=\max_i\{\tilde{\sigma}^2_i\}$. 
Therefore, from Lemma \ref{lem:var_bd} we conclude that:
\begin{align}\label{eq:var_bd:snet}
\mathbb{E}_{\xi,\phi}[\|{ \bar{G}_{\mu,i}(x^r,\phi^r,\xi^r)}-\nabla f_{\mu,i}(x^r)\|^2]\le \frac{\tilde{\sigma}^2}{J}.
\end{align}

Let us define the auxiliary sequence $y^r:=\{y_i^r\}_{i=1}^N$ as follows

	\begin{align}\label{eq:def:y}
	y^0=x^0, \; y_j^r = y_j^{r-1}, \; \mbox{if} \; j\ne i_r, \; \mbox{else} \; \; y^r_{i_r}=x^r,~ \forall~r. 
	\end{align}
	Next let us define the potential function which measures the progress of algorithm 
	
\begin{align*}
\tilde{Q}^r &= \sum_{i=1}^{N} f_{\mu,i}(x^r)+\sum_{i=1}^{N} \frac{4}{\alpha_i\rho_i}\|\nabla f_{\mu,i}(y_i^{r-1}) - \nabla f_{\mu,i}(x^r) \|^2
 + h(x^r),
\end{align*} where $f_{\mu,i}(x^r)$ denotes the smoothed version of function $f_i(x^r)$ { as defined similarly} in \eqref{eq:grad_smooth_ver}.

First, we study the behavior of the potential function.  For this algorithm let us define the filtration $\cF^{r}$ as the $\sigma$-field generated by $\{i_t,\phi^t,\xi^t\}_{t=1}^{r-1}$. Throughout this section the expectations are taken with respect to { $\{i_{r},\phi^{r},\xi^{r}\}$} conditioning on $\cF^{r}$ unless otherwise noted. 
\begin{lemma}\label{lem:zone-ns:des}
	Suppose  Assumptions C holds true. Set $\tilde{p}:=\sum_{i=1}^{N}\frac{1}{p_i}$, $\beta:=\frac{1}{\sum_{i=1}^N \rho_i}$, and for $i=1,2,\cdots, N$, we pick
	
	\begin{align}\label{eq:p:eta:zo}
	\alpha_i= p_i =\frac{\rho_i}{\sum_{i=1}^N{\rho_i}}, ~ \mbox{and} ~ \rho_i\geq \frac{5.5L_{\mu,i}}{p_i}, ~  i=1,\cdots N.
	\end{align}
	 Then we have the following result for the ZONE-S algorithm
		\begin{align}\label{eq:zone-ns:descent}
		\mathbb{E}[\tilde{Q}^{r+1}-\tilde{Q}^{r}]&\le \frac{-1}{100\beta}\mathbb{E}\|x^{r+1}-x^{r}\|^2\nonumber\\ 
		&- \sum_{i=1}^N\frac{1}{2\rho_i}\mathbb{E}\|\nabla f_{\mu,i}(x^{r})-\nabla f_{\mu,i}(y_i^{r-1})\|^2  +\frac{3\tilde{p}\beta\tilde{\sigma}^2}{J}.
		\end{align}
\end{lemma}
Next we define the optimality gap as the following
	\begin{align}\label{eq:gap:def:zo}
	\Psi^r :=\frac{1}{\beta^2}\mathbb{E}\bigg\|x^r-\prox_h^{1/\beta}[x^r-\beta \nabla f(x^r) ]\bigg\|^2,
	\end{align}
where ${\rm \prox}^{\gamma}_{h}[u]:=\argmin\;\;h(u)+\frac{\gamma}{2}\|x-u\|^2$ is the proximity operator for function $h$. Note that when the nonsmooth term $h\equiv0$, $\Psi^r$ reduces to the size of the gradient vector $\mathbb{E}\|\nabla f(x^r)\|^2$.
\begin{remark}
	From the parameter selection in \eqref{eq:p:eta:zo} one can derive the following relationships (see \cite[Theorem 2.1]{davood_NIPS2016_nestt}):
	\begin{align}
	p_i=\frac{\sqrt{L_{\mu,i}}}{\sum_{i=1}^N\sqrt{L_{\mu,i}}},\; \; \rho_i=\sqrt{5.5L_{\mu,i}}\sum_{i=1}^{N}\sqrt{5.5L_{\mu,i}}.
	\end{align}
	In particular, the probability of picking agent $i_r$ is not uniform. Utilizing this nonuniform sampling we are able to improve the algorithm speed. See \cite[Sec. 4]{davood_NIPS2016_nestt} for detailed discussion.  
\end{remark}
Finally we present the main convergence result about the proposed ZONE-S algorithm. 
\begin{theorem}\label{thm:rate:zonestt-g}
	Suppose Assumptions C hold true, and $u$ is uniformly randomly sampled from $\{1,2,\cdots, T\}$. Let us set
	 $\frac{1}{\beta} = {5.5(\sum_{i=1}^N\sqrt{L_{\mu,i}})^2}$. Then we have the following bounds for the optimality gap in expectation
		\begin{align*}
		&1)~ \mathbb{E}_{u}[\Psi^u]\leq {  2200}\Big(\sum_{i=1}^N\sqrt{L_{\mu,i}}\Big)^2\frac{\mathbb{E}[{\tilde{Q}^1-\tilde{Q}^{T+1}} ]}{T} 
		+ \frac{\mu^2L^2(M+3)^3}{2}+ \frac{1024\tilde{p}\tilde{\sigma}}{J};\\
		&2)~\mathbb{E}_{u}[\Psi^u] + \mathbb{E}_u\bigg[\sum_{i=1}^N{3 \rho^2_i }\bigg\|z^{u}_i-x^{u-1}\bigg\|^2\bigg]
		\leq {  2200} \bigg(\sum_{i=1}^N\sqrt{L_{\mu,i}}\bigg)^2\frac{\mathbb{E}[{\tilde{Q}^1-\tilde{Q}^{T+1}} ]}{T}\nonumber\\
		&+\frac{\mu^2L^2(M+3)^3}{2}+ \frac{1024\tilde{p}\tilde{\sigma}}{J}.
		\end{align*}
\end{theorem}  Note that part (1) only measures the primal optimality gap, while part (2) also shows that the expected constraint violation shrinks in the same order.

\begin{remark}
	Similar to the ZONE-M, the bound for the optimality gap of ZONE-S is dependent on two $T$-independent constants, the first one $\frac{\mu^2{L}^2(M+3)^3}{2}$  arises from using zeroth-order gradient, and the second term $\frac{1024\tilde{p}\tilde{\sigma}^2}{J}$ arise from the uncertainty in the gradient estimation.  Again, if we pick $\mu\in\mathcal{O}(\frac{1}{\sqrt{T}})$, and $J\in\mathcal{O}({T})$, we obtain the following sublinear convergence rate
	\begin{align}
	\mathbb{E}_{u}[\Psi^u]&\leq {  2200}\Big(\sum_{i=1}^N\sqrt{{L}_{\mu,i}}\Big)^2\frac{\mathbb{E}[{\tilde{Q}^1-\tilde{Q}^{T+1}} ]}{T} +\frac{1024\tilde{p}\tilde{\sigma}^2}{T}
	+ \frac{{L}^2(M+3)^3}{T}.
	\end{align}
\end{remark}
\begin{remark}
	{The reason that the ZONE-S is able to incorporate non-smooth terms, in contrast to the ZONE-M algorithm,  is that it has special network structure. In particular, the non-smooth term is optimized by the central controller, and the fact that the central controller can talk to every node makes sure that the non-smooth term is optimized by using the most up-to-date information from the network. }
\end{remark}

\section{Numerical Results}\label{sec:numerical}
In this section we numerically evaluate the effectiveness of the ZONE-M and ZONE-S algorithms.   We consider some distributed nonconvex optimization problems in zeroth-order setup (i.e., we only have access to the noisy functional values).  We set the noise $\xi$ to be a  zero-mean Gaussian random  variable with standard deviation $\sigma = 0.01$. All the simulations are performed on Matlab 2015a on a Laptop with 4 GB memory and Intel Core i7-4510U CPU (2.00 GHz), running on Linux (Ubuntu 16.04) operating system. 
\subsection{ZONE-M Algorithm}
\begin{figure}
	\centering
	\begin{subfigure}{.4\textwidth}
		\hspace{.3cm}\includegraphics[width=2.2 in]{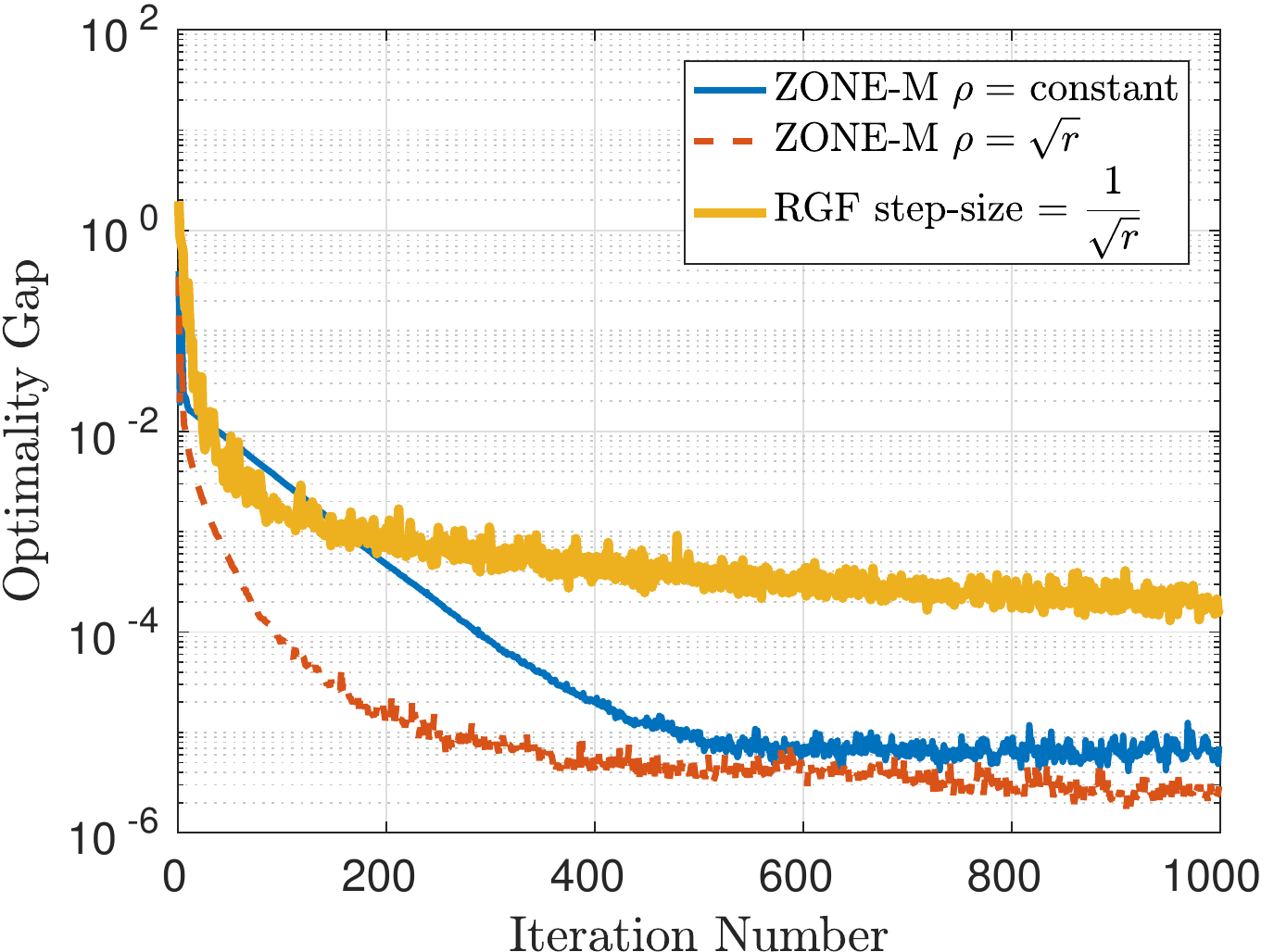}
		\caption{\footnotesize The $x$-axis represents the  number of algorithm iterations, and the $y$-axis measures the optimality gap defined in \eqref{eq:crit}. }
		\label{fig:opt_gap_conc} 
	\end{subfigure}
	\hfill
	\begin{subfigure}{.4\textwidth}
		\hspace{.3cm}\includegraphics[width=2.25 in]{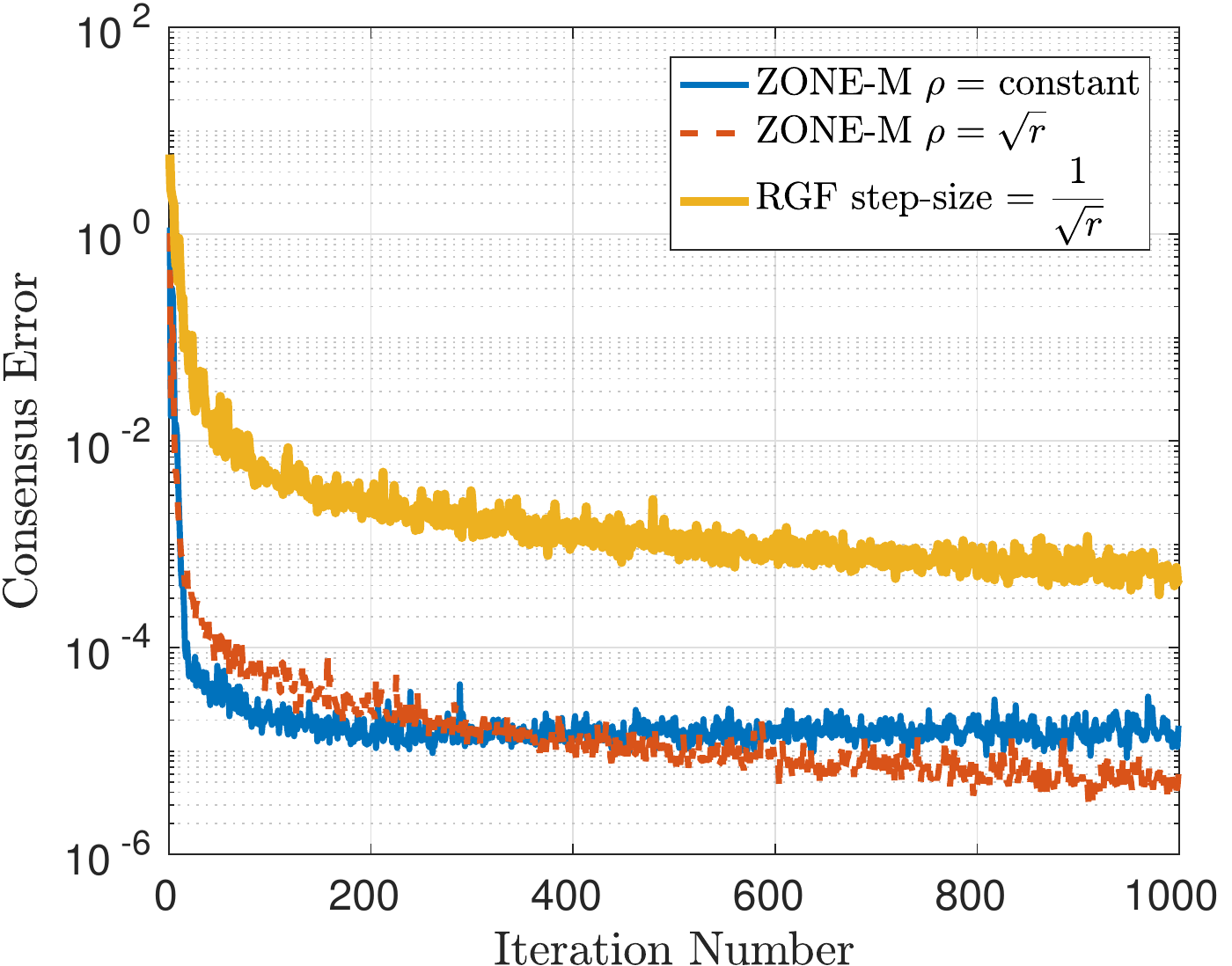}
		\caption{\footnotesize The $x$-axis represents the  number of algorithm iterations, and the $y$-axis measures the constraint violation defined in \eqref{eq:crit}.}
		\label{fig:conc_vio_conc} 
	\end{subfigure}   
	\caption{\footnotesize Comparison of different algorithms for the problem \eqref{p1_problem}.}
	\label{fig:concensus}   
\end{figure}

We study the following nonconvex distributed optimization problems. 
Consider minimizing sum of nonconvex functions in a distributed setting
\begin{equation} \label{p1_problem}
\mathop {\min }\limits_{z\in \mathbb{R}^Q} \sum\limits_{i = 1}^N {{f_i}({z_i})} ,\quad \textrm{s.t.}\;Az = 0.
\end{equation}
where each agent $i$ can only obtain the zeroth-order information of its local function, given by {$${f_i}({z_i}) = \frac{{a_i}}{{1 + {e^{ - {z_i}}}}} + b_i\textrm{log}(1 + z_i^2),$$ }
where $a_i$ and $b_i$ are constants generated from an i.i.d Gaussian distribution. Clearly the function $f_i$ is nonconvex and smooth, {  and we can simply check that it satisfies assumption A, B}. In our experiments the graphs are generated based on the scheme proposed in  \cite{YildizScag08}. In this scheme  a random graph with $N$ nodes and radius $R$ is generated  with nodes uniformly distributed over a unit square, and two nodes connect if their distance is less than $R$. We set problem dimension {$M=1$}, and the number of nodes in the network $N=20$ with radius $R=0.6$.  {  The penalty parameter $\rho$ is selected to satisfy theoretical bounds given in Lemma \ref{lem:bd:pot}}, the smoothing parameter is set $\mu=\frac{1}{\sqrt{T}}$, and we set $J=T$, where maximum number of iterations is picked $T=1000$.
{We compare the ZONE-M algorithm with Randomized Gradient Free (RGF) algorithm  with diminishing stepsize $\frac{1}{\sqrt{r}}$ ($r$ denotes the iterations counter) proposed in \cite{yuan2015randomized}, which is only developed for convex problems.  
	We also compare our algorithm with a variant of ZONE-M which uses increasing penalty parameter $\rho=\sqrt{r}$. When choosing $\rho=\sqrt{r}$ neither RGF not ZONE-M has convergence guarantee}.
We use the optimality gap (opt-gap) and constraint violation (cons-vio), displayed below, to measure the quality of the solution generated by different algorithms
\begin{align}\label{eq:crit}
\mbox{opt-gap}&:={\bigg\| {\sum\limits_{i = 1}^N {{\nabla}{f_i}({z_i})} } \bigg\|^2} + {\big\| A{z} \big\|^2},\nonumber\\
\mbox{cons-vio}&:=\|A{  z}\|^2.
\end{align}
Figure \ref{fig:concensus} illustrates the comparison among different algorithms. Each point in the figure is obtained by averaging over 50 independent trials. One can observe that: 1) ZONE-M  converges faster compared with RGF in both the optimality gap and the consensus error; 2) ZONE-M with increasing penalty ($\rho = \sqrt{r}$) appears to be faster than its constant stepsize counterpart.  

In the next set of experiments we compare different algorithms with a number of choices of network size, i.e., $N\in \{10,20,40,80\}$. For this problem we set the radius $R=0.5$ 
The results (average over 50 independent trials) are reported in Table \ref{tab:cons}. In this table { ZONE-M (C)} and {ZONE-M (I)} denote ZONE-M with constant and increasing penalty parameter, respectively.  
We observe that ZONE-M algorithm is always faster compared with the RGF. 
\begin{table}[]
	\centering
	\caption{\footnotesize Comparison results for ZONE-M and RGF}
	\label{tab:cons}
	\begin{tabular}{lcccccc}
		\cline{2-7}
		& \multicolumn{3}{c|}{\textbf{opt-gap}}                & \multicolumn{3}{c}{\textbf{cons-error}} \\ \hline
		N  & ZONE-M(C) & ZONE-M(I)          & \multicolumn{1}{l|}{RGF} & ZONE-M(C)  & ZONE-M(I)           & RGF    \\ \hline
		10  & 6.8E-6  & \textbf{8.8E-6} & 1.7E-4                   & 2.5E-5   & \textbf{2.0E-5}  & 0.002  \\
		20 & 4.2E-5  & \textbf{2.2E-5} & 5.3E-3                    & 3.1E-5   & \textbf{2.2E-5}  & 0.003  \\
		40 & 7.0E-5  & \textbf{3.0E-5} & 1.8E-3                    & {3.8E-4}   & \textbf{2.8E-4}  & 0.017  \\
		80 & 5.7E-4  & \textbf{7.5E-5} & 0.014                      & {5.4E-4}   & \textbf{3.0E-4}    & 0.09  \\ \hline      
	\end{tabular}
\end{table}

\subsection{ZONE-S Algorithm}
In this subsection we { demonstrate the performance }of the ZONE-S algorithm. The penalty parameter $\rho$ is selected to satisfy the conditions given in Lemma \ref{lem:zone-ns:des}, or to be an increasing sequence satisfying $\rho=\sqrt{r}$. 
For comparison purpose we consider two additional algorithms, namely the zeroth-order gradient descent (ZO-GD)  \cite{nesterov2011random} ({ which is a centralized algorithm}), and the zeroth-order stochastic gradient descent (ZO-SGD) \cite{ghadimi2013stochastic}. To be notationally consistent with our algorithm we denote the stepsize for these two algorithms with $1/\rho$. For ZO-GD  it has been shown that {if the stepsize is set $1/\rho=\frac{1}{4L(M+4)}$, and the smoothing factor satisfies $\mu\le\mathcal{O}(\frac{\epsilon}{ML})$, then the algorithm will converge to an $\epsilon$-stationary solution \cite[Section 7]{nesterov2011random}}. {Also, for ZO-SGD the optimality gap decreases in the order of $\frac{1}{\sqrt{T}}$ when we pick stepsize $1/\rho<\frac{1}{2(M+4)}$, and the smoothing parameter $\mu$ satisfies $\mu\le \frac{D_f}{(M+4)\sqrt{2N}}$, where $D_f:=\bigg[\frac{2(f(x^1)-f^*)}{L}\bigg]^{1/2}$ ($f^*$ denotes the optimal value)} \cite[Theorem 3.2]{ghadimi2013stochastic}.  Note that the theoretical results for ZO-SGD is valid only for smooth cases, however we include it here for comparison purposes.

\noindent{\bf Nonconvex Sparse Optimization Problem.}
Consider the following optimization problem 
\begin{align}
&\min_{x\in \mathbb{R}^M} \;\sum_{i=1}^{N} f_i(x) ~\st~ \|x\|_1\leq \ell,
\end{align}
where $f_i(x)=x^\top \Gamma x -\gamma^\top x$, ($\Gamma \in \mathbb{R}^{M\times M}$, and $\gamma\in \mathbb{R}^M$), and $\ell$ is a positive constant that controls the sparsity level of the solution. In this problem the matrix $\Gamma\in\mathbb{R}^{M\times M}$ is not necessarily a positive semidefinite matrix, thus the problem is not convex; see for example high dimensional regression problem with noisy observations in \cite[problem (2.4)]{loh12} . This problem is a special case of the original problem in \eqref{eq:sum:nonsmoth} with $h(x)$ being the indicator function of the set $\{x \mid \|x\|_1\leq \ell\}$. 

\begin{figure}[t!]
	\centering
	\includegraphics[width=2.2 in]{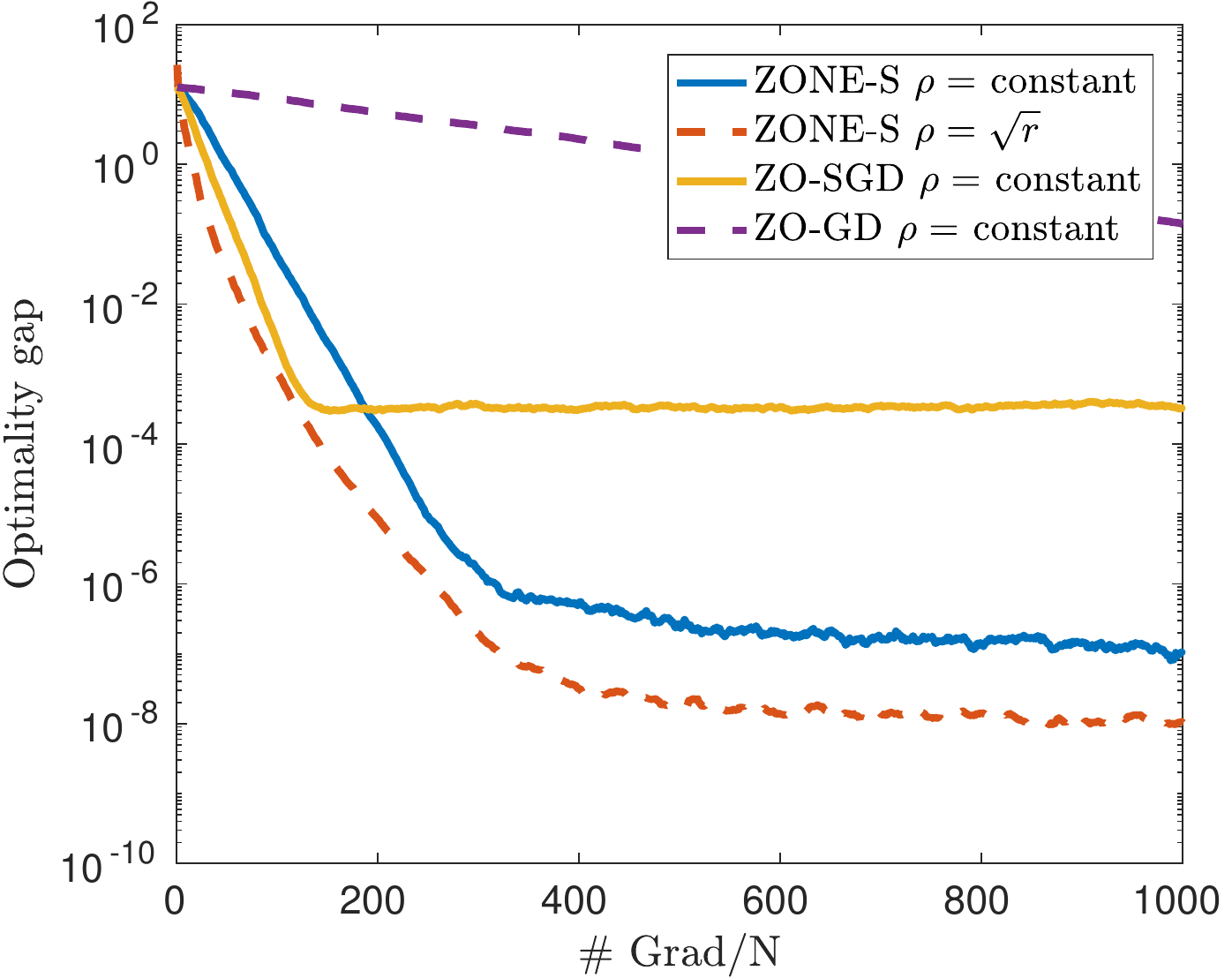}
	\caption{\footnotesize  The Optimality Gap for Nonconvex Sparse Optimization problem. The $x$-axis represents the effective
		passes through the data, which is the number of gradient evaluations divided by the  number of agents in the network. The $y$-axis measures the optimality gap defined in \eqref{eq:gap:def:zo}. }
	\label{fig:opt_gap:reg_l1}
\end{figure}

We compare the following four algorithms: ZONE-S with constant stepsize  $$\rho_i = \sqrt{5.5L_{\mu,i}}\sum_{i=1}^{N}\sqrt{5.5L_{\mu,i}};$$ ZONE-S with increasing penalty parameter $\rho_i=\sqrt{r}$; ZO-GD with  constant stepsize ($1/\rho=\frac{1}{4L(M+4)}$), and ZO-SGD with  constant step size $1/\rho=\frac{1}{2L(M+4)}$.  The problem dimension is set as $N=10$, and $M=100$.  The algorithm stops when the iteration counter reaches $T=1000$. The results are plotted  in Figure \ref{fig:opt_gap:reg_l1}, which depicts the progress of the optimality gap [defined as in \eqref{eq:gap:def:zo}] versus the number of iterations. Each point in this figure is obtained by averaging over 50 independent trials. We can observe that ZONE-S converges faster than the ZO-GD and ZO-SGD. Furthermore, the performance of ZONE-S improves when using the increasing stepsize, as compared to that of the constant stepsize.  

\section{Conclusion} \label{sec:con}
In this work, we consider nonconvex multi-agent  optimization problem under zeroth-order setup. We design algorithms to solve the problem over two popular network structures, namely MNet and SNet. 
We have rigorously analyzed the convergence rate of the proposed algorithms and we have proved that both algorithms converge to the set of first-order stationary solutions under very mild conditions on the problem and by appropriately choosing the algorithm parameters.


\section*{APPENDIX A} 
In this appendix we provide the proofs for ZONE-M. 
\noindent\subsection{\bf Proof of Lemma \ref{lemma:mu:bound}}
 Rearranging terms in  \eqref{eq:lambda:zero} we get 
	\begin{align}\label{eq:due:rearg}
	\lambda^{r+1}-\lambda^r=\rho Az^{r+1}.
	\end{align}
 Utilizing this equation and equation \eqref{eq:kkt:z:compact}, we obtain
\begin{align}\label{eq:kkt:z:app}
G_{\mu}^{J,r}+A^\top\lambda^{r+1}+\rho L^{+}(z^{r+1}-z^r)=0.
\end{align}
From \eqref{eq:due:rearg} it is clear that $\lambda^{r+1}-\lambda^r$ lies in the column space of $A$, therefore the following is true
\begin{align}\label{eq:bd:sig}
\sqrt{\sigma_{\min}}\|\lambda^{r+1}-\lambda^r\|\leq\|A^\top (\lambda^{r+1}-\lambda^r)\|,
\end{align} 
where $\sigma_{\min}$ denotes the smallest non-zero eigenvalue of $A^\top A$. 

Replacing $r$ with $r-1$ in equation \eqref{eq:kkt:z:app}, and then using the definition of $w^r:=(z^{r+1}-z^r) -(z^{r}-z^{r-1})$ we obtain
	\begin{align}\label{eq:lam}
	&\bigg\|\lambda^{r+1}-\lambda^{r}\bigg\|^2\le \frac{1}{{\sigma_{\min}}}\bigg\|G_{\mu}^{J,r}-G_{\mu}^{J,r-1}{ +}\rho L^{+}w^r\bigg\|^2\nonumber\\
	&=\frac{1}{{\sigma_{\min}}}\|G_{\mu}^{J,r}-G_{\mu}^{J,r-1}+\nabla g_{\mu}(z^r)-\nabla g_{\mu}(z^r){ +}\rho L^{+}w^r\|^2\nonumber\\
	&\stackrel{\eqref{eq:rel2}}\le \frac{3}{{\sigma_{\min}}}\|G_{\mu}^{J,r}-\nabla g_{\mu}(z^r)\|^2+\frac{3}{{\sigma_{\min}}}\|\nabla g_{\mu}(z^r)- G_{\mu}^{J,r-1}\|^2
	+\frac{3\rho^2}{\sigma_{\min}}\| L^{+}w^r\|^2.
	\end{align}
	 Let us add and subtract $\nabla g_{\mu}(z^{r-1})$ to the second term on the RHS of \eqref{eq:lam}, and 
take the expectation on both sides ( expectations are taken with respect to $(\xi^{r+1}, \phi^{r+1})$ conditioning on filtration $\mathcal{F}^{r}$ defined previously.)
\begin{align}
&\mathbb{E}\|\lambda^{r+1}-\lambda^{r}\|^2\le \frac{3}{{\sigma_{\min}}}\mathbb{E}\|G^{J,r}_{\mu}-\nabla g_{\mu}(z^r)\|^2 \nonumber\\
&+\frac{6}{\sigma_{\min}}\mathbb{E}\|\nabla g_\mu(z^r)-\nabla g_\mu(z^{r-1})\|^2+\frac{3\rho^2}{\sigma_{\min}}\mathbb{E}\| L^{+}w^r\|^2
\nonumber\\
&+\frac{6}{{\sigma_{\min}}}\mathbb{E}\|\nabla g_{\mu}(z^{r-1})- G_{\mu}^{J,r-1}\|^2\nonumber\\
&\stackrel{\rm(i)}\le \frac{9\tilde{\sigma}_g^2}{J\sigma_{\min}}+\frac{6L^2_\mu}{{\sigma_{\min}}}\mathbb{E}	\|z^r-z^{r-1}\|^2+\frac{3\rho^2\|L^{+}\|}{\sigma_{\min}}\mathbb{E}\|w^r\|_{L^+}^2,	
\end{align} where $\rm(i)$ is true because of \eqref{eq:var_bd:mnet}, and the facts that $\nabla g_\mu(z)$ is  ${L}_\mu$-smooth  and $\|L^+w^r\|^2\le \|L^+\|\|w^r\|^2_{L^{+}}$. The lemma is proved. \QED
\noindent\subsection{\bf Proof of Lemma \ref{lem:bd:pot}}
{ Using Assumption B.1, and the  fact that $D\succeq I$, it can be shown that if $2\rho\ge \hat{L}$, then function
	\begin{align*}
	&L_\rho(z,\lambda)+\frac{\rho}{2}\|z-z^r\|^2_{L^{+}}
	=g(z)+\langle \lambda,Az\rangle+\frac{\rho}{2}\|Az\|^2+\frac{\rho}{2}\|z-z^r\|^2_{L^{+}},
	\end{align*} 
	  is strongly convex  with modulus $2\rho-\hat{L}$.   See \cite[Theorem 2.1]{Zlobec2005}.
Using this fact, let us  bound  $L_\rho^{r+1}-L_\rho^r$. 
	\begin{align}
	&L_{\rho}^{r+1} - L_{\rho}^r= L_{\rho}^{r+1} - L_{\rho}(z^{r+1},\lambda^{r})+L_{\rho}(z^{r+1},\lambda^{r})- L_{\rho}^r\nonumber\\
	&\stackrel{\rm(i)}\leq \langle\nabla_z L_{\rho}(z^{r+1}, \lambda^r) + \rho L^{+}(z^{r+1}-z^r),z^{r+1}-z^{r}\rangle\nonumber\\
	&+\frac{1}{\rho}\|\lambda^{r+1}-\lambda^r\|^2-\frac{2\rho-\hat{L}}{2}\|z^{r+1}-z^r\|^2 . \label{eq:lag:bound:derive}
	\end{align} where $\rm(i)$ is true due to the  strong convexity of $L_\rho(z,\lambda)+\frac{\rho}{2}\|z-z^r\|^2_{L^{+}}$  with modulus $2\rho-\hat{L}$ and \eqref{eq:due:rearg}. Now using \eqref{eq:kkt:z:app} we further have
	\begin{align}
	&L_{\rho}^{r+1} - L_{\rho}^r\le \big\langle \nabla g(z^{r+1})-G_{\mu}^{J,r},z^{r+1}-z^r\big\rangle\nonumber\\
	&+\frac{1}{\rho}\|\lambda^{r+1}-\lambda^r\|^2-\frac{2\rho-\hat{L}}{2}\|z^{r+1}-z^r\|^2\nonumber\\
	&\stackrel{\rm(i)}\le \frac{1}{\rho}\|\lambda^{r+1}-\lambda^r\|^2+\frac{{ \hat{L}^2}-2\rho+\hat{L}}{2}\|z^{r+1}-z^r\|^2\nonumber\\
	&+\frac{1}{2\hat{L}^2}\|\nabla g(z^{r+1})-G_{\mu}^{J,r}\|^2\nonumber,
	\end{align} where $\rm(i)$ is application of \eqref{eq:rel3} for $\epsilon=\hat{L}^2$. Taking expectation on both sides we get	

	\begin{align}
	&\mathbb{E}\bigg[L_{\rho}^{r+1} - L_{\rho}^r\bigg]\stackrel{\rm(i)}\le \frac{9\tilde{\sigma}_g^2}{\rho J\sigma_{\min}}+\frac{6L^2_\mu}{{\rho\sigma_{\min}}}\mathbb{E}\|z^r-z^{r-1}\|^2\nonumber\\
	&+\frac{3\rho\|L^+\|}{\sigma_{\min}}\mathbb{E}\|w^r\|_{L^+}^2+\frac{\hat{L}^2-2\rho+\hat{L}}{2}\mathbb{E}\|z^{r+1}-z^r\|^2
	+\frac{1}{2\hat{L}^2}\mathbb{E}\|\nabla g(z^{r+1})-G_{\mu}^{J,r}\|^2\nonumber\\
   &\stackrel{\rm(ii)}\le\bigg(\frac{9}{\rho\sigma_{\min}}+\frac{3}{2\hat{L}^2}\bigg)\frac{\tilde{\sigma}_g^2}{J}+\frac{3\mu^2(Q+3)^3}{8}+\frac{6L^2_\mu}{{\rho\sigma_{\min}}}\mathbb{E}\|z^r-z^{r-1}\|^2\nonumber\\
   &+\frac{3\rho\|L^+\|}{\sigma_{\min}}\mathbb{E}\|w^r\|_{L^+}^2+\frac{\hat{L}^2-2\rho+\hat{L}+3}{2}\mathbb{E}\|z^{r+1}-z^r\|^2, \label{eq:aug:diff}
	\end{align} where in $\rm(i)$ we use Lemma \ref{lemma:mu:bound} to bound $\mathbb{E}\|\lambda^{r+1}-\lambda^r\|$,  in $\rm(ii)$ we apply \eqref{eq:rel2}, \eqref{eq:diff:grad:bd}, \eqref{eq:var_bd:mnet} , and the fact that $\nabla g_\mu(z)$ is $L_\mu$-smooth with $L_\mu\le\hat{L}$.
	
	Next we bound $V^{r+1}-V^r$.  Applying  the optimality condition for problem \eqref{eq:z:zero} together with equation \eqref{eq:lambda:zero}  yields the following
	\begin{align*}
	\langle &G_{\mu}^{J,r}+A^\top \lambda^{r+1} + \rho L^{+}(z^{r+1}-z^r),z^{r+1}-z\rangle\le 0,\; \forall~z\in\mathbb{R}^Q.
	\end{align*}
	Similarly, for the $(r-1)$th iteration, we have 
	\begin{align*}
	\langle &G_{\mu}^{J,r-1}+A^\top \lambda^{r}+ \rho L^{+}(z^{r}-z^{r-1}),z^{r}-z\rangle\le 0,\; \forall~z\in\mathbb{R}^Q.
	\end{align*}
	Now let us set $z=z^r$ in first, $z=z^{r+1}$ in second equation, and add them. We obtain
	\begin{align}\label{eq:sum:kkt}
	&\langle A^\top (\lambda^{r+1}-\lambda^r), z^{r+1}-z^r\rangle\le
	-\langle G_{\mu}^{J,r} - G_{\mu}^{J,r-1}+ \rho L^{+}w^r, z^{r+1}-z^r\rangle.
	\end{align}
   The left hand side can be expressed in the following way
	\begin{align}\label{eq:lhs}
	&\langle  A^\top (\lambda^{r+1}-\lambda^r), z^{r+1}-z^r\rangle=\rho\langle  Az^{r+1}, A z^{r+1}-Az^r\rangle\nonumber\\
	&\stackrel{\eqref{eq:rel1}}=\frac{\rho}{2}\bigg(\|Az^{r+1}\|^2-\|Az^{r}\|^2+\|A(z^{r+1}-z^r)\|^2\bigg).
	\end{align}
For the right hand side we have

	\begin{align*}
	&-\langle G_{\mu}^{J,r} - G_{\mu}^{J,r-1}+ \rho L^{+} w^r,z^{r+1}-z^r\rangle\nonumber\\
	 &=-\langle G_{\mu}^{J,r} - G_{\mu}^{J,r-1}, z^{r+1}-z^r\rangle - \langle\rho L^{+} w^r, z^{r+1}-z^r\rangle \nonumber\\
	&\stackrel{\eqref{eq:rel3}}\le  \frac{1}{2\hat{L}}\|G_{\mu}^{J,r} - G_{\mu}^{J,r-1}\|^2+\frac{\hat{L}}{2}\|z^{r+1}-z^r\|^2
	-\rho \langle L^{+} w^r, z^{r+1}-z^r\rangle\\
	&\stackrel{\rm(i)}\le \frac{3}{2\hat{L}}\bigg(\|G^{J,r}_\mu-\nabla g_\mu(z^{r})\|^2+\|\nabla g_\mu(z^{r-1})-G^{J,r-1}_\mu\|^2\nonumber\\
	&+\|\nabla g_\mu(z^{r})-\nabla g_\mu(z^{r-1})\|^2\bigg)+\frac{\hat{L}}{2}\|z^{r+1}-z^r\|^2 
	-\rho \langle L^{+} w^r, z^{r+1}-z^r\rangle.
	\end{align*}
	To get $\rm(i)$  we add and subtract $\nabla g_\mu(z^{r}) + \nabla g_\mu(z^{r-1})$ to  $G_{\mu}^{J,r} - G_{\mu}^{J,r-1}$ and use \eqref{eq:rel2}.   
	Taking expectation on both sides, we have

	\begin{align}\label{eq:rhs}
	&-\mathbb{E}[\langle G_{\mu}^{J,r} - G_{\mu}^{J,r-1}+ \rho L^{+}w^r, z^{r+1}-z^r\rangle]\nonumber\\
	&\stackrel{\eqref{eq:var_bd:mnet}}\le \frac{3}{2\hat{L}}\bigg(\frac{2\tilde{\sigma}_g^2}{J}+\hat{L}^2\mathbb{E}\|z^r-z^{r-1}\|^2\bigg)+\frac{\hat{L}}{2}\mathbb{E}\|z^{r+1}-z^r\|^2
	-\rho \mathbb{E}[\langle L^{+} w^r, z^{r+1}-z^r\rangle]\nonumber\\	
	&\stackrel{\rm(i)}=\frac{3\tilde{\sigma}_g^2}{\hat{L}J}+\frac{3\hat{L}}{2}\mathbb{E}\|z^{r} - z^{r-1}\|^2+\frac{\hat{L}}{2}\mathbb{E}\|z^{r+1}-z^r\|^2\nonumber\\
	&+\frac{\rho}{2}\mathbb{E}\bigg[\|z^{r}-z^{r-1}\|^2_{L^{+}}-\|z^{r+1}-z^r\|^2_{L^{+}} -\|w^r\|_{L^{+}}^2\bigg],
	\end{align}
	 where in $\rm(i)$ we apply \eqref{eq:rel1} with $b=(L^+)^{1/2}(z^{r+1}-z^r)$ and $a=(L^+)^{1/2}(z^{r}-z^{r-1})$.  Combining \eqref{eq:sum:kkt}, \eqref{eq:lhs} and \eqref{eq:rhs}, we obtain 
	 \begin{align}\label{eq:rl}
	 &\frac{\rho}{2}\mathbb{E}\bigg(\|Az^{r+1}\|^2-\|Az^{r}\|^2+\|A(z^{r+1}-z^r)\|^2\bigg)\nonumber\\
	 &
	 \le \frac{3\tilde{\sigma}_g^2}{\hat{L}J}+\frac{3\hat{L}}{2}\mathbb{E}\|z^{r} - z^{r-1}\|^2+\frac{\hat{L}}{2}\mathbb{E}\|z^{r+1}-z^r\|^2\nonumber\\
	 &+\frac{\rho}{2}\mathbb{E}\bigg(\|z^{r}-z^{r-1}\|^2_{L^{+}}-\|z^{r+1}-z^r\|^2_{L^{+}} -\|w^r\|_{L^{+}}^2\bigg).
	 \end{align}
	 Recall that matrix $B:=L^++\frac{k}{c\rho}I_Q$, and $V^{r+1}$ is defined as
	 	\begin{align*}
	 	V^{r+1}&:=\frac{\rho}{2}\bigg(\|Az^{r+1}\|^2+\|z^{r+1}-z^r\|^2_{B}\bigg)\nonumber\\
	 	&=\frac{\rho}{2}\bigg(\|Az^{r+1}\|^2+\|z^{r+1}-z^r\|^2_{L^+}\bigg) +\frac{k}{2c}\|z^{r+1}-z^r\|^2.
	 	\end{align*}
	Rearranging terms in \eqref{eq:rl}, we have
	\begin{align}\label{eq:vpart}
	&\mathbb{E}[V^{r+1}-V^r]\le \bigg(\frac{\hat{L}}{2}+\frac{k}{2c}\bigg)\mathbb{E}\|z^{r+1}-z^r\|^2+\frac{3\tilde{\sigma}_g^2}{\hat{L}J}\nonumber\\
	& + \bigg(\frac{3\hat{L}}{2}-\frac{k}{2c}\bigg)\mathbb{E}\|z^{r} - z^{r-1}\|^2
	-\frac{\rho}{2}\mathbb{E}\bigg(\|w^r\|_{L^{+}}^2+\|A(z^{r+1}-z^r)\|^2\bigg)\nonumber\\
	&\le \bigg(\frac{\hat{L}}{2}+\frac{k}{2c}\bigg)\mathbb{E}\|z^{r+1}-z^r\|^2+ \bigg(\frac{3\hat{L}}{2}-\frac{k}{2c}\bigg)\mathbb{E}\|z^{r} - z^{r-1}\|^2
-\frac{\rho}{2}\mathbb{E}\|w^r\|_{L^{+}}^2 +\frac{3\tilde{\sigma}_g^2}{\hat{L}J}.
	\end{align}
Now let us consider the definition of $P^{r+1}:=L^{r+1}_{\rho}+ cV^{r+1}$. 
	Utilizing \eqref{eq:aug:diff}, and \eqref{eq:vpart} and definition of $k$ as  $	 k:=2\bigg(\frac{6\hat{L}^2}{\rho\sigma_{\min}}+\frac{3c\hat{L}}{2}\bigg)$ eventually we obtain	
	\begin{align}
&\mathbb{E}\bigg[P^{r+1}-P^{r}\bigg]\le \frac{k-c_1}{2}\mathbb{E}\|z^{r+1}-z^r\|^2
-c_2\mathbb{E}\|w^r\|_{L^{+}}^2+\frac{3\mu^2({Q}+3)^3}{8} +c_3\frac{\tilde{\sigma}_g^2}{J},
\end{align} 
where we define,
\begin{align*}
&c_1 := 2\rho-\hat{L}^2-(c+1)\hat{L}-3,\nonumber\\
&c_2:=\bigg(\frac{c \rho}{2}-\frac{3\rho\|L^{+}\|}{\sigma_{\min}}\bigg), \quad c_3:=\frac{9}{\rho\sigma_{\min}}+\frac{3}{2\hat{L}^2}+\frac{3c}{\hat{L}}.
\end{align*}
The lemma is proved.
\QED
\noindent\subsection{\bf Proof of Lemma \ref{lemma:lower:bound}}
From \eqref{eq:kkt:z:app} we have
 \begin{align}\label{eq:opt_z}
 G_{\mu}^{J,r}+A^\top\lambda^{r+1}+\rho L^{+}(z^{r+1}-z^r)=0.
 \end{align}
{From this equation we have
 \begin{align}
 \|A^\top\lambda^{r+1}\|^2&= \|G_{\mu}^{J,r}+\rho L^{+}(z^{r+1}-z^r)\|^2
\leq 2\|G_{\mu}^{J,r}\|^2+2\rho^2 \|L^{+}(z^{r+1}-z^r)\|^2.
 \end{align}
  Also from the fact that $\lambda^0=0$
we have that the dual variable  lies in the column space of $A$ and one can conclude that 
 \begin{align}
 \sigma_{\min}\|\lambda^{r+1}\|^2\leq 2\|G_{\mu}^{J,r}\|^2+2\rho^2 \|L^{+}(z^{r+1}-z^r)\|^2.
 \end{align}}
 Dividing both sides by $\sigma_{\min}$ yields
 \begin{align}\label{lam}
 \|\lambda^{r+1}\|^2\leq \frac{2}{\sigma_{\min}}\|G_{\mu}^{J,r}\|^2+\frac{2\rho^2}{\sigma_{\min}}\|L^+(z^{r+1}-z^r)\|^2.
 \end{align}
 Now based on the definition of potential function we have
 \begin{align}\label{pot}
 P^{r+1} &= g(z^{r+1}) +\frac{\rho}{2}\|Az^{r+1}+\frac{1}{\rho}\lambda^{r+1}\|^2-\frac{1}{2\rho}\|\lambda^{r+1}\|^2\nonumber\\
 &+\frac{c\rho}{2}\|Az^{r+1}\|^2 +\frac{c\rho}{2}\|z^{r+1}-z^r\|_B^2,
 \end{align}
 where {$B:=L^++\frac{k}{c\rho}I$} [note that $k=2(\frac{6\hat{L}^2}{\rho\sigma_{\min}}+\frac{3c\hat{L}}{2})$]. Plugging \eqref{lam} in \eqref{pot}, and utilizing the fact that $g(z^{r+1})\geq 0$ from Assumption [B2], $\frac{c\rho}{2}\|Az^{r+1}\|^2\geq 0$, and $\|Az^{r+1}+\frac{1}{\rho}\lambda^{r+1}\|^2\geq 0$  we get
 \begin{align*}
 P^{r+1} &\geq \frac{-1}{\rho\sigma_{\min}}\|G_{\mu}^{J,r}\|^2-\frac{\rho}{\sigma_{\min}}\|L^{+}(z^{r+1}-z^r)\|^2
 +\frac{c\rho}{2}\|z^{r+1}-z^r\|_B^2\\
 &\stackrel{\rm(i)}\geq \frac{-1}{\rho\sigma_{\min}}\|G_{\mu}^{J,r}\|^2 + \frac{\rho}{\sigma_{\min}}\|z^{r+1}-z^r\|^2_{\frac{c\sigma_{\min}}{2}L^+-(L^+)^2},
 \end{align*}
 where $\rm(i)$ is true because $\frac{k}{2}\|z^{r+1}-z^r\|^2\geq0$. Notice that $L^+$ is a symmetric PSD matrix. Therefore, picking constant $c$ large enough such that 
 	$c\geq \frac{2\|L^+\|}{\sigma_{\min}}$, 
 we have $\frac{c\sigma_{\min}}{2}L^+-(L^+)^2\succeq 0$. Hence, with this choice of $c$ we get the following bound for the potential function
 \begin{align}
 P^{r+1} \geq -\frac{1}{\rho\sigma_{\min}}\|G_{\mu}^{J,r}\|^2.
 \end{align}
 Taking expectation on both sides we have
 \begin{align}\label{pot2}
 \mathbb{E}[P^{r+1}] \geq -\frac{1}{\rho\sigma_{\min}}\mathbb{E}\|G_{\mu}^{J,r}\|^2.
 \end{align}
 Now let us prove that $\mathbb{E}\|G_{\mu}^{J,r}\|^2$ is upper bounded as follows:
 \begin{align}
 \mathbb{E}\|G_{\mu}^{J,r}\|^2 &= \mathbb{E}\|G_{\mu}^{J,r}-\nabla g_\mu(z^r)+\nabla g_\mu(z^r)\|^2\nonumber\\
 &\leq 2\mathbb{E}\|G_{\mu}^{J,r}-\nabla g_\mu(z^r)\|^2 + 2\mathbb{E}\|\nabla g_\mu(z^r)\|^2\nonumber\\
 &\stackrel{\rm(i)}\leq \frac{2\tilde{\sigma}^2_g}{J} + 2\mathbb{E}\|\nabla g_\mu(z^r)\|^2\nonumber\\
 &\stackrel{\rm(ii)}\leq 2\tilde{\sigma}^2_g + 4\mathbb{E}\|\nabla g(z^r)\|^2+\mu^2\hat{L}^2(Q+3)^3\nonumber\\
 &\stackrel{\rm(iii)} \leq 2\tilde{\sigma}^2_g + 4K^2+\mu^2\hat{L}^2(Q+3)^3
 \end{align}
 where $\rm(i)$ is true due to \eqref{eq:var_bd:mnet}, $\rm(ii)$ comes from the fact that $J\geq 1$, and $\|\nabla g_\mu(z^r)\|^2\leq 2\|\nabla g(z^r)\|^2+\frac{\mu^2}{2}\hat{L}^2(Q+3)^3$ \cite[Theorem 3.1]{ghadimi2013stochastic}, and in $\rm(iii)$ we use assumption A1 in the paper in which we assumed there exists a $K$ such that $\|\nabla g(z)\|\leq K$. Therefore, we have proved that there exists a constant $K_2:=2\tilde{\sigma}^2_g + 4K^2+\mu^2\hat{L}^2(Q+3)^3$ such that $\mathbb{E}\|G_{\mu}^{J,r}\|^2\leq K_2$. Finally, plugging this bound in equation \eqref{pot2}, we get 
 \begin{align}\label{pot3}
 \mathbb{E}[P^{r+1}] \geq -\frac{1}{\rho\sigma_{\min}}K_2. 
 \end{align}
 Since $K_2$ is not dependent on $T$, in order to prove the Lemma we just need to set $T$-independent lower bound $\underline{P}:=-\frac{1}{\rho\sigma_{\min}}K_2$. 

\noindent\subsection{\bf Proof of Theorem \ref{thm:conv}}
 Let us bound the optimality gap given in \eqref{eq:opt_gap} term by term.  First we bound the gradient of AL function with respect to variable $z$ in point $(z^{r+1}, \lambda^r)$ in the following way
		\begin{align}
		&\|\nabla_z L_\rho(z^{r+1}, \lambda^r)\|^2 = \|\nabla g(z^{r+1})+A^\top \lambda^r +\rho A^\top Az^{r+1}\|^2\nonumber\\
		&\stackrel{\eqref{eq:lambda:zero}}=\|\nabla g(z^{r+1})+A^\top \lambda^{r+1}\|^2\nonumber\\
		&\stackrel{\eqref{eq:opt_z}}=\|\nabla g(z^{r+1})-G_{\mu}^{J,r}-\rho L^{+}(z^{r+1}-z^r)\|^2\nonumber\\
		&\stackrel{\eqref{eq:rel2}}\le 2\|\nabla g(z^{r+1})-G_{\mu}^{J,r}\|^2+2\rho^2\|L^{+}(z^{r+1}-z^r)\|^2\nonumber\\
		&\stackrel{\rm(i)}\le 4(\|\nabla g(z^{r+1}) -\nabla g_{\mu}(z^r)\|^2+\|\nabla g_{\mu}(z^{r})-G_{\mu}^{J,r}\|^2)
		+2\rho^2\|L^{+}(z^{r+1}-z^r)\|^2,\label{eq:bd:grad1}
		\end{align} where in $\rm(i)$ we add and subtract $\nabla g_\mu(z^r)$ to $\nabla g(z^{r+1})-G_{\mu}^{J,r}$ and apply \eqref{eq:rel2} and \eqref{eq:kkt:z:app}.
	Further, let us take expectation on both sides of \eqref{eq:bd:grad1}
		\begin{align}
		&\mathbb{E}\|\nabla_z L_\rho(z^{r+1}, \lambda^r)\|^2 \nonumber\\
		&\le 4\mathbb{E}\bigg(\|\nabla g(z^{r+1}) -\nabla g_{\mu}(z^r)\|^2+\|\nabla g_\mu(z^{r})-G^{J,r}_{\mu}\|^2\bigg)
		+2\rho^2\mathbb{E}\|L^{+}(z^{r+1}-z^r)\|^2\nonumber\\
		&\stackrel{\rm(i)}\le 8\mathbb{E}\bigg(\|\nabla g(z^{r+1}) -\nabla g_{\mu}(z^{r+1})\|^2+\hat{L}^2\|z^{r+1}-z^r\|^2\bigg)+\frac{4\tilde{\sigma}^2}{J}
		+2\rho^2\mathbb{E}\|L^{+}(z^{r+1}-z^r)\|^2\nonumber\\
		&\stackrel{\eqref{eq:diff:grad:bd}}\le 2\mu^2\hat{L}^2(Q+3)^3 +8\hat{L}^2\mathbb{E}\|z^{r+1}-z^r\|^2+\frac{4\tilde{\sigma}_g^2}{J}
		+2\rho^2\mathbb{E}\|L^{+}(z^{r+1}-z^r)\|^2, \label{eq:bd:grad2}
		\end{align} where in  $\rm(i)$ we applied \eqref{eq:var_bd:mnet}, \eqref{eq:rel2}, and the fact that  $\nabla g_\mu(z)$ is $L_\mu$-smooth with $L_\mu\le\hat{L}$.
	Second, let us bound the expected value of the constraint violation. Utilizing the equation  \eqref{eq:lambda:zero} we have
		\begin{align*}
		\|Az^{r+1}\|^2=\frac{1}{\rho^2}\|\lambda^{r+1}-\lambda^r\|^2.
		\end{align*}
	Taking expectation  on the above identity, and utilizing the fact that $L_\mu\le \hat{L}$, and \eqref{eq:mu:difference:bound}, we obtain the following
		\begin{align}\label{eq:bd:vio}
		\mathbb{E}\|Az^{r+1}\|^2&=\frac{1}{\rho^2}\mathbb{E}\|\lambda^{r+1}-\lambda^r\|^2\le\frac{9\tilde{\sigma}_g^2}{J\rho^2\sigma_{\min}}\nonumber\\
		&+\frac{6\hat{L}^2}{\rho^2\sigma_{\min}}\mathbb{E}\|z^{r}-z^{r-1}\|^2+\frac{3\|L^{+}\|}{\sigma_{\min}}\mathbb{E}\|w^r\|_{L^{+}}^2.
		\end{align}
	Summing up  \eqref{eq:bd:grad2} and \eqref{eq:bd:vio}, we have the following bound for the optimality gap 
		\begin{align}
		\Phi^{r+1}&\le \alpha_1\mathbb{E}\|z^{r+1}-z^r\|^2+\alpha_2\mathbb{E}\|z^r-z^{r-1}\|^2+\alpha_3\mathbb{E}\|w^r\|_{L^{+}}^2\nonumber\\
		&+\bigg(\frac{9+4\rho^2\sigma_{\min}}{\rho^2\sigma_{\min}}\bigg)\frac{\tilde{\sigma}_g^2}{J}+2\mu^2\hat{L}^2({Q}+3)^3,\label{eq:Q:bd}
		\end{align}
	where $\alpha_1, \alpha_2, \alpha_3$ are positive constants given by 
		$$\alpha_1=8\hat{L}^2+2\rho^2{  \|L^{+}\|^2}, \; \alpha_2=\frac{6\hat{L}^2}{\rho^2\sigma_{min}}, \; \alpha_3=\frac{3\|L^{+}\|}{\sigma_{min}}.$$ 
	Summing both sides of \eqref{eq:Q:bd},  we obtain the following
		\begin{align}
		\sum_{r=1}^{T}\Phi^{r+1}&\le \sum_{r=1}^{T-1}(\alpha_1+\alpha_2)\mathbb{E}\|z^{r+1}-z^{r}\|^2+\sum_{r=1}^{T} \alpha_3\mathbb{E}\| w^r\|_{L^{+}}^2\nonumber\\
		&+\alpha_2\mathbb{E}\|z^1-z^0\|^2 +\alpha_1\mathbb{E}\|z^{{T}+1}-z^{T}\|^2\nonumber\\
		&+2T\mu^2\hat{L}^2(Q+3)^3+T\bigg(\frac{9+4\rho^2\sigma_{\min}}{\rho^2\sigma_{\min}}\bigg)\frac{\tilde{\sigma}_g^2}{J}.\label{eq:sum:Q}
		\end{align}
	Applying Lemma \ref{lem:bd:pot} and  summing both sides of \eqref{eq:diff:bd}  over $T$ iterations, we obtain
		\begin{align}
		&\mathbb{E}\bigg[P^1-P^{T+1}\bigg]\ge \sum_{r=1}^{T-1}\frac{c_1-k}{2}\mathbb{E}\|z^{r+1}-z^{r}\|^2+\sum_{r=1}^{T} c_2\mathbb{E}\| w^r\|_{L^+}^2\nonumber\\
		&+\frac{c_1-k}{2}\mathbb{E}\|z^{T+1}-z^T\|^2-\frac{3T\mu^2(Q+3)^3}{8}-\frac{Tc_3\tilde{\sigma}_g^2}{J}.\label{eq:sum:P}
		\end{align}
	Let us set $\zeta=\frac{\max(\alpha_1+\alpha_2, \alpha_3)}{\min(\frac{c_1-k}{2}, c_2)}$. Combining the two inequalities \eqref{eq:sum:Q} and \eqref{eq:sum:P}, and utilizing the fact that  $\mathbb{E}[P^{T+1}]$ is lower bounded by $\underline{P}$, we arrive at the following inequality
		\begin{align}\label{eq:sum:Q:2}
		&\sum_{r=1}^{T}\Phi^{r+1}\le  \zeta\mathbb{E}[P^1-\underline{P}]+\alpha_2\mathbb{E}\|z^1-z^0\|^2\nonumber\\
		&+T\bigg(\zeta c_3+\frac{9+4\rho^2\sigma_{\min}}{\rho^2\sigma_{\min}}\bigg)\frac{\tilde{\sigma}_g^2}{J}+T\bigg( \frac{3\zeta}{8}+2\hat{L}^2\bigg)\mu^2{  (Q+3)^3}.
		\end{align}
	Since $u$ is a uniformly random variable in the set $\{1,2,\cdots, T\}$ we have
		\begin{align}\label{eq:exp:J}
		\mathbb{E}_u[\Phi^u]=\frac{1}{T}\sum_{r=1}^{T}\Phi^{r+1}.
		\end{align}
	Dividing both sides of \eqref{eq:sum:Q:2} on $T$ and using  \eqref{eq:exp:J} implies the following
		\begin{align*}
		\mathbb{E}_u[\Phi^u]&\le \frac{\zeta\mathbb{E}[P^1-\underline{P}]+\alpha_2\mathbb{E}\|z^1-z^0\|^2}{T}\nonumber\\
		&+ \bigg(\zeta c_3+\frac{9+4\rho^2\sigma_{\min}}{\rho^2\sigma_{\min}}\bigg)\frac{\tilde{\sigma}_g^2}{J}+\bigg(\frac{3\zeta}{8}+2\hat{L}^2\bigg)\mu^2{  (Q+3)^3}
		\end{align*}
	By setting 
		\begin{align}
		&\gamma_1=\zeta\mathbb{E}[P^1-\underline{P}]+\alpha_2\mathbb{E}\|z^1-z^0\|^2,\nonumber\\
		& \gamma_2=\zeta c_3+\frac{9+4\rho^2\sigma_{\min}}{\rho^2\sigma_{\min}},~ \gamma_3 =\bigg(\frac{3\zeta}{8}+2\hat{L}^2\bigg){  (Q+3)^3},
		\end{align}
	we conclude the proof. \QED

\section*{APPENDIX B} 
This appendix contains the proof of the lemmas in Section \ref{sec:snet} which are related ZONE-S. 

In order to facilitate the derivations, in the following let us present some key properties of ZONE-S. Let us define $r(j) := \max\{t\mid t<r+1, j = i_t\}$ which is the { most recent} iteration in which  agent $j$ is picked before iteration $r+1$. From this definition we can see that $r(i_r)=r$. Let us repeat the update equations of ZONE-S algorithm
\begin{align}
z_{i_r}^{r+1}&=x^r-\frac{1}{\alpha_{i_r}\rho_{i_r}}\bigg[\lambda^r_{i_r}+\bar{G}_{\mu,{i_r}}(x^r,\phi^r,\xi^r)\bigg];\label{eq:z_i:zone-ns:app}&&\\
\lambda^{r+1}_{i_r}&=\lambda_{i_r}^{r}+\alpha_{i_r}\rho_{i_r}\bigg(z^{r+1}_{i_r}-x^{r}\bigg); \label{eq:lam:zone-ns:app}&&\\
\lambda^{r+1}_{j}&=\lambda_{j}^{r}, \quad z^{r+1}_{j} = x^r, \quad \forall~j\ne {i_r}. \label{eq:z_i:zone-ns2:app}&&
\end{align}
{\bf \noindent Property 1: Compact form for dual update}. 
Combining \eqref{eq:z_i:zone-ns:app}, \eqref{eq:lam:zone-ns:app}, and using the definition of $r(j)$  we get
\begin{align}
&\lambda_{i_r}^{r+1}=-\bar{G}_{\mu,i_r}(x^r,\phi^r,\xi^r),\\
&\lambda_j^{r+1}=\lambda_j^r=-\bar{G}_{\mu,j}(x^{r(j)},\phi^{r(j)},\xi^{r(j)}), \quad j\neq i_r.
\end{align}
Using the definition of sequence $y^r$ [$y^0=x^0, \; y_j^r = y_j^{r-1}, \; \; \mbox{if} \; \; j\ne i_r, \; \mbox{else} \; \; y^r_{i_r}=x^r, ~ \forall~r \ge 1$] we have $y_i^r = x^{r(i)}$ for all $i=1,2,\cdots N$. Using this we get the following compact form
\begin{align}
\lambda^{r+1}_{i} & = -\bar{G}_{\mu,i}(y^r_i, \phi^{r(i)},\xi^{r(i)}), \; \forall i\label{eq:lambda:compact}.
\end{align}
{\bf \noindent Property 2: Compact form for primal update}. 
From \eqref{eq:z_i:zone-ns2:app}, and \eqref{eq:lambda:compact} for $j\neq i_r$ we have
\begin{align}
z_j^{r+1}=x^r\stackrel{\eqref{eq:lambda:compact}}=& x^r-\frac{1}{\alpha_j\rho_j}[\lambda_j^{r+1}+\bar{G}_{\mu,j}(y^r_j, \phi^{r(j)},\xi^{r(j)})]\nonumber\\
\stackrel{\eqref{eq:z_i:zone-ns2:app}}=& x^r-\frac{1}{\alpha_j\rho_j}[\lambda_j^{r}+\bar{G}_{\mu,i}(y^r_j, \phi^{r(j)},\xi^{r(j)})].\label{eq:z_j}
\end{align}
Considering \eqref{eq:z_i:zone-ns:app}, and \eqref{eq:z_j} we can express the update equation for $z$ in ZONE-S algorithms in the following compact form 
\begin{align}
z^{r+1}_{i} & = x^r -\frac{1}{\alpha_{i}\rho_{i}} \bigg[\lambda^r_{i} +\bar{G}_{\mu,i}(y^r_i, \phi^{r(i)},\xi^{r(i)})\bigg] ,\forall i\label{eq:x:compact}.
\end{align}
{\bf \noindent Property 3: Bound the distance between update direction and the gradient direction}. Let us define 
\begin{align}
u^{r+1} &:= \beta\bigg(\sum_{i=1}^{N} \rho_i z^{r+1}_i +\sum_{i=1}^{N}{\lambda^{r}_i}\bigg), \label{eq:def:u}
\end{align}}where we set $\beta:= 1/\sum_{i=1}^{N}\rho_i$. Using \eqref{eq:def:u}, it is easy to check that $x$-update \eqref{eq:x:zone-ns} is equivalent to solving the following problem
\begin{align}\label{eq:x:prox}
x^{r+1} &= \arg\min_{x}\; \frac{1}{2\beta}\|x-u^{r+1}\|^2 + h(x)\nonumber\\
&=\prox_h^{1/\beta}(u^{r+1}).
\end{align}
The optimality condition for this problem is given by
\begin{align}
x^{r+1}-u^{r+1}+\beta \eta^{r+1} =0, \label{eq:opt:x:zone-ns}
\end{align} where $\eta^{r+1}\in \partial h(x^{r+1})$ is a subgradient of  {  $h$ at $x^{r+1}$}.    [When there is no confusion we use the shorthand notation $\bar{G}_{\mu,i}^r$ to denote $\bar{G}_{\mu,i}(x^r,\phi^r,\xi^r)$]
\begin{align}\label{eq:u2}
&u^{r+1} = \beta\bigg(\sum_{i=1}^{N} \rho_i z^{r+1}_i +\sum_{i=1}^{N}{\lambda^{r}_i}\bigg)\nonumber\\
&\stackrel{\eqref{eq:z_i:zone-ns2:app}}=\beta\bigg(\sum_{i=1}^{N}\rho_ix^r-\rho_{i_r}(x^r-z_{i_r}^{r+1})+\sum_{i=1}^{N}\lambda_i^r\bigg)\nonumber\\
&\stackrel{\eqref{eq:lambda:compact}, \eqref{eq:z_i:zone-ns:app}}=x^r-\frac{\beta}{\alpha_{i_r}}\bigg[\bar{G}_{\mu,i_r}^r-\bar{G}_{\mu,i_r}(y_{i_r}^{r-1},\phi^{(r-1)(i_r)},\xi^{(r-1)(i_r)})\bigg]\nonumber\\
&\qquad-\beta\sum_{i=1}^{N}\bar{G}_{\mu,i}(y_i^{r-1},\phi^{(r-1)(i)},\xi^{(r-1)(i)}).
\end{align}
Let us  further define 
\begin{align}\label{eq:v_ir}
&v_{i_{r}}^{r} := \sum_{i=1}^N \bar{G}_{\mu,i}(y_i^{r-1}, \phi^{(r-1)(i)},\xi^{(r-1)(i)})\nonumber\\
&\qquad+\frac{1}{ \alpha_{i_{r}}}\bigg[\bar{G}_{\mu,i_{r}}^r- \bar{G}_{\mu,i_{r}}(y^{r-1}_{i_{r}}, \phi^{(r-1)(i_r)},\xi^{(r-1)(i_r)})\bigg].
\end{align}
We conclude that 
\begin{align}
u^{r+1}=x^r - \beta v^{r}_{i_{r}}\label{eq:u}.
\end{align}
Plugging \eqref{eq:u} into \eqref{eq:opt:x:zone-ns} we obtain 
\begin{align}
x^{r+1}=x^r-\beta(v^{r}_{i_{r}} +\eta^{r+1})\label{eq:x:ex}.
\end{align}
From the  definition \eqref{eq:v_ir} it is clear that $v_{i_r}^r$ is an approximation of certain gradient of $\sum_{i=1}^{N} f_{i}(x^r)$. Below we make this intuition precise by bounding the  $\|\sum_{i=1}^{N}\nabla f_{\mu,i}(x^r)-v_{i_r}^r\|$.   Using the definition of $v_{i_r}^r$ we have
\begin{align}
&\bigg\|\sum_{i=1}^{N}\nabla f_{\mu,i}(x^r)-v_{i_r}^r\bigg\|^2\nonumber\\
&\stackrel{  \eqref{eq:v_ir}}=\bigg\|\sum_{i=1}^{N}\nabla f_{\mu,i}(x^r) - \sum_{i=1}^N \bar{G}_{\mu,i}(y_i^{r-1}, \phi^{(r-1)(i)},\xi^{(r-1)(i)})\nonumber\\
&\quad-\frac{1}{ \alpha_{i_{r}}}\bigg[\bar{G}_{\mu,i_{r}}^r- \bar{G}_{\mu,i_{r}}(y^{r-1}_{i_{r}}, \phi^{(r-1)(i_r)},\xi^{(r-1)(i_r)})\bigg]\bigg\|^2.
\end{align}
Let us set $\cJ^{r}:=\{i_r,\phi^r,\xi^r\}$. Setting $\alpha_{i}=p_i$ and taking conditional expectation on both sides, we have 
\begin{align}
&\mathbb{E}_{\cJ^r}\bigg[\bigg\|\sum_{i=1}^{N}\nabla f_{\mu,i}(x^r)-v_{i_r}^r\bigg\|^2|\mathcal{F}^r\bigg]\nonumber\\
&= \mathbb{E}_{\cJ^r}\bigg[\bigg\|\sum_{i=1}^{N}\nabla f_{\mu,i}(x^r) - \sum_{i=1}^N \bar{G}_{\mu,i}(y_i^{r-1}, \phi^{(r-1)(i)},\xi^{(r-1)(i)})\nonumber\\
&\quad-\frac{1}{ \alpha_{i_{r}}}\bigg[\bar{G}_{\mu,i_{r}}^r- \bar{G}_{\mu,i_{r}}(y^{r-1}_{i_{r}}, \phi^{(r-1)(i_r)},\xi^{(r-1)(i_r)})\bigg]\bigg\|^2\mathcal{F}^r\bigg]\nonumber\\
&\stackrel{\rm(i)}\le \mathbb{E}_{\cJ^r}\bigg[\bigg\|\frac{\bar{G}_{\mu,i_{r}}^r- \bar{G}_{\mu,i_{r}}(y^{r-1}_{i_{r}}, \phi^{(r-1)(i_r)},\xi^{(r-1)(i_r)})}{\alpha_{i_{r}}}\bigg\|^2\mathcal{F}^r\bigg], \nonumber
\end{align}  where $\rm(i)$ is true because  $\mathbb{E}[\|x-\mathbb{E}[x]\|^2]=\mathbb{E}{[\|x\|^2]}-\|\mathbb{E}[x]\|^2\leq\mathbb{E}{[\|x\|^2]}$ and the following identity
\begin{align*}
&\mathbb{E}_{\cJ^r}\bigg[\frac{1}{ \alpha_{i_{r}}}\bigg[\bar{G}_{\mu,i_{r}}^r- \bar{G}_{\mu,i_{r}}(y^{r-1}_{i_{r}}, \phi^{(r-1)(i_r)},\xi^{(r-1)(i_r)})\bigg]\mid\cF^{r}\bigg]\\
&=\sum_{i=1}^{N}\nabla f_{\mu,i}(x^r) - \sum_{i=1}^N \bar{G}_{\mu,i}(y_i^{r-1}, \phi^{(r-1)(i)},\xi^{(r-1)(i)}).
\end{align*}
Now if we take expectation with respect to $i_r$, (given $\cF^{r}$)
\begin{align*}
&\mathbb{E}_{\cJ^r}\bigg[\bigg\|\sum_{i=1}^{N}\nabla f_{\mu,i}(x^r)-v_{i_r}^r\bigg\|^2|\mathcal{F}^r\bigg]\nonumber\\
&\leq\sum_{i=1}^{N}\frac{1}{p_i}\mathbb{E}_{\phi^r,\xi^r}\bigg[\bigg\|\bar{G}_{\mu,i}^r- \bar{G}_{\mu,i}(y^{r-1}_{i}, \phi^{(r-1)(i)},\xi^{(r-1)(i)})\bigg\|^2|\mathcal{F}^r\bigg] \nonumber\\
&= \sum_{i=1}^{N}\frac{1}{p_i}\mathbb{E}_{\phi^r,\xi^r}\bigg[\bigg\|\bar{G}_{\mu,i}^r-\nabla f_{\mu,i}(x^r)+\nabla f_{\mu,i}(x^r)-\nabla f_{\mu,i}(y_i^{r-1})\nonumber\\
&\quad+\nabla f_{\mu,i}(y_i^{r-1})-\bar{G}_{\mu,i}(y^{r-1}_{i}, \phi^{(r-1)(i)},\xi^{(r-1)(i)})\bigg\|^2|\mathcal{F}^r\bigg].
\end{align*}
Then utilizing \eqref{eq:rel2} and \eqref{eq:var_bd:snet}, we have
\begin{align}
&\mathbb{E}_{\cJ^r}\bigg[\bigg\|\sum_{i=1}^{N}\nabla f_{\mu,i}(x^r)-v_{i_r}^r\bigg\|^2|\mathcal{F}^r\bigg]\nonumber\\
&\le 3\sum_{i=1}^{N}\frac{1}{p_i}\bigg(\frac{\tilde{\sigma}^2}{J}+\bigg\|\nabla f_{\mu,i}(x^r)-\nabla f_{\mu,i}(y_i^{r-1})\bigg\|^2\nonumber\\
&\quad+\bigg\|\nabla f_{\mu,i}(y_i^{r-1})- \bar{G}_{\mu,i}(y^{r-1}_{i}, \phi^{(r-1)(i)},\xi^{(r-1)(i)})\bigg\|^2\bigg).
\end{align}
 Using the definition of  $\tilde{p}=\sum_{i=1}^{N}\frac{1}{p_i}$, overall we have the following
\begin{align}\label{eq:vir_bd:cond}
&\mathbb{E}_{\cJ^r}\bigg[\bigg\|\sum_{i=1}^{N}\nabla f_{\mu,i}(x^r)-v_{i_r}^r\bigg\|^2\mid\cF^{r}\bigg]\nonumber\\
&\le  \frac{3\tilde{p}\tilde{\sigma}^2}{J}+\sum_{i=1}^{N}\frac{3}{p_i}\bigg(\bigg\|\nabla f_{\mu,i}(x^r)-\nabla f_{\mu,i}(y_i^{r-1})\bigg\|^2\nonumber\\
&\quad+\bigg\|\nabla f_{\mu,i}(y_i^{r-1})- \bar{G}_{\mu,i}(y^{r-1}_{i}, \phi^{(r-1)(i)},\xi^{(r-1)(i)})\bigg\|^2\bigg) .
\end{align}
Using the property of conditional expectation we have 
\begin{align}\label{eq:cond:pro}
&\mathbb{E}\bigg\|\sum_{i=1}^{N}\nabla f_{\mu,i}(x^r)-v_{i_r}^r\bigg\|^2=\mathbb{E}_{\cF^{r},\cJ^r}\bigg\|\sum_{i=1}^{N}\nabla f_{\mu,i}(x^r)-v_{i_r}^r\bigg\|^2\nonumber\\
&=  \mathbb{E}_{\cF^r}\bigg[\mathbb{E}_{\cJ^r}\big[\big\|\sum_{i=1}^{N}\nabla f_{\mu,i}(x^r)-v_{i_r}^r\big\|^2\mid\cF^{r}\big]\bigg].
\end{align}
Now let us break the filtration as $\cF^r=\cF^r _1\cup\cF^r _2$ where   $\cF^r _1:=\{i_t\}_{t=1}^{r-1}$, and $\cF^r _2:=\{\phi^t,\xi^t\}_{t=1}^{r-1}$. Using these notations we have
\begin{align}\label{eq:Q:r-1}
&\mathbb{E}_{\cF^r}\bigg\|\nabla f_{\mu,i}(y_i^{r-1})- \bar{G}_{\mu,i}(y^{r-1}_{i}, \phi^{(r-1)(i)},\xi^{(r-1)(i)})\bigg\|^2\nonumber\\
&= \mathbb{E}_{\cF_1^r}\bigg[\mathbb{E}_{\cF_2^r}\bigg\|\nabla f_{\mu,i}(y_i^{r-1})-\bar{G}_{\mu,i}(y^{r-1}_{i}, \phi^{(r-1)(i)},\xi^{(r-1)(i)})\bigg\|^2\mid \cF^r_1\bigg]\nonumber\\
&\stackrel{\eqref{eq:var_bd:snet}}\le \frac{\tilde{\sigma}^2}{J}.
\end{align}
Combining \eqref{eq:vir_bd:cond}, \eqref{eq:cond:pro}, \eqref{eq:Q:r-1}, we obtain
\begin{align}\label{eq:vir_bd}
&\mathbb{E}\bigg\|\sum_{i=1}^{N}\nabla f_{\mu,i}(x^r)-v_{i_r}^r\bigg\|^2\nonumber\\
&\le \sum_{i=1}^{N}\frac{3}{p_i}\mathbb{E}\bigg\|\nabla f_{\mu,i}(x^r)-\nabla f_{\mu,i}(y_i^{r-1})\bigg\|^2+\frac{6\tilde{p}\tilde{\sigma}^2}{J}.
\end{align}
\noindent\subsection{\bf Proof of Lemma \ref{lem:zone-ns:des}}
By assumption $\alpha_i=p_i$, according to the definition of potential function $\tilde{Q}^r$, we have
\begin{align}
&\mathbb{E}_{\cJ^r}[\tilde{Q}^{r+1}- \tilde{Q}^{r}\mid \cF^{r}]\nonumber\\
&=\mathbb{E}_{\cJ^r}\bigg[\sum_{i=1}^N\bigg(f_{\mu,i}(x^{r+1})-f_{\mu,i}(x^{r})\bigg) +h(x^{r+1})-h(x^{r})\mid {\cF^{r}}\bigg]\nonumber\\
&+\mathbb{E}_{\cJ^r}\bigg[\sum_{i=1}^N\frac{4}{p_i\rho_i}\bigg\|\nabla f_{\mu,i}(x^{r+1})-\nabla f_{\mu,i}(y_i^{r})\bigg\|^2-\frac{4}{p_i\rho_i}\bigg\|\nabla f_{\mu,i}(x^{r})-\nabla f_{\mu,i}(y_i^{r-1})\bigg\|^2\mid  {\cF^{r}}\bigg]\label{eq:p:diff}.
\end{align}
The proof consists of the following steps:\\
{\bf Step 1).}  We  bound the first term in \eqref{eq:p:diff} as follows 
\begin{align}\label{eq:diff:p:1}
&\mathbb{E}_{\cJ^r}\bigg[\sum_{i=1}^N \bigg(f_{\mu,i}(x^{r+1})-f_{\mu,i}(x^{r})\bigg)+h(x^{r+1})-h(x^{r})\mid {\cF^{r}}\bigg]\nonumber\\
&\stackrel{\rm (i)}\leq \mathbb{E}_{\cJ^r}\bigg[\sum_{i=1}^N\langle \nabla f_{\mu,i}(x^{r}),x^{r+1}-x^{r}\rangle+\langle \eta^{r+1}, x^{r+1}-x^{r}\rangle+\frac{\sum_{i=1}^N{L_{\mu,i}}}{2}\|x^{r+1}-x^{r}\|^2\mid {\cF^{r}}\bigg]\nonumber\\
&=\mathbb{E}_{\cJ^r}\bigg[\big\langle \sum_{i=1}^N \nabla f_{\mu,i}(x^{r})+\eta^{r+1} +\frac{1}{\beta} (x^{r+1}-x^{r}),x^{r+1}-x^{r}\big\rangle\mid {\cF^{r}}\bigg]\nonumber\\
&-\bigg(\frac{1}{\beta}-\frac{\sum_{i=1}^{N}L_{\mu,i}}{2}\bigg)\mathbb{E}_{\cJ^r}\bigg[\|x^{r+1}-x^{r}\|^2\mid {\cF^{r}}\bigg],\nonumber
\end{align}
 where in $\rm(i)$ we have used the Lipschitz continuity of $\nabla f_{\mu,i}$ as well as the convexity of $h$.  Then from \eqref{eq:x:ex} we further have
\begin{align}
&\mathbb{E}_{\cJ^r}\bigg[\sum_{i=1}^N \bigg(f_{\mu,i}(x^{r+1})-f_{\mu,i}(x^{r})\bigg)+h(x^{r+1})-h(x^{r})\mid {\cF^{r}}\bigg]\nonumber\\
&\le\mathbb{E}_{\cJ^r}\bigg[\big\langle\sum_{i=1}^N \nabla f_{\mu,i}(x^{r})-v_{i_{r}}^{r},x^{r+1}-x^{r}\big\rangle\mid {\cF^{r}}\bigg]\nonumber\\
&\quad-\bigg(\frac{1}{\beta}-\frac{\sum_{i=1}^{N}L_{\mu,i}}{2}\bigg)\mathbb{E}_{\cJ^r}\bigg[\|x^{r+1}-x^{r}\|^2\mid {\cF^{r}}\bigg]\nonumber\\
&\stackrel{\rm(i)}\le \sum_{i=1}^{N}\frac{3\beta}{2p_i}\bigg(\bigg\|\nabla f_{\mu,i}(x^r)-\nabla f_{\mu,i}(y_i^{r-1})\bigg\|^2\nonumber\\
&\quad+\bigg\|\nabla f_{\mu,i}(y_i^{r-1})- \bar{G}_{\mu,i}(y^{r-1}_{i}, \phi^{(r-1)(i)},\xi^{(r-1)(i)})\bigg\|^2\bigg) \nonumber\\
&-\bigg(\frac{1}{2\beta}-\frac{\sum_{i=1}^NL_{\mu,i}}{2}\bigg)\mathbb{E}_{\cJ^r}\bigg[\|x^{r+1}-x^{r}\|^2\mid {\cF^{r}}\bigg] +\frac{3\tilde{p}\beta\tilde{\sigma}^2}{2J},
\end{align} where  in $\rm(i)$  we utilize \eqref{eq:rel3} with $\epsilon=\frac{1}{\beta}$, and  \eqref{eq:vir_bd:cond}. 

	{\bf Step 2).} In this step we bound the second term in equation \eqref{eq:p:diff} as follows 
	\begin{align}
	&\mathbb{E}_{\cJ^r}\bigg[\big\|\nabla f_{\mu,i}(x^{r+1})-\nabla f_{\mu,i}(y_i^{r})\big\|^2\mid {\cF^{r}}\bigg]\nonumber\\
	&\stackrel{\rm (i)}\leq (1+\epsilon_i)\mathbb{E}_{\cJ^r}\bigg[\big\|\nabla f_{\mu,i}(x^{r+1})-\nabla f_{\mu,i}(x^{r})\big\|^2\mid {\cF^{r}}\bigg]\nonumber\\
	&+\bigg(1+\frac{1}{\epsilon_i}\bigg)\mathbb{E}_{\cJ^r}\bigg\|\nabla f_{\mu,i}(x^{r})-\nabla f_{\mu,i}(y_i^{r})\|^2\mid {\cF^{r}}\bigg]\label{eq:del:fmui0}\\
	&\stackrel{\rm (ii)}= (1+\epsilon_i)\mathbb{E}_{\cJ^r}\bigg[\|\nabla f_{\mu,i}(x^{r+1})-\nabla f_{\mu,i}(x^{r})\|{ ^2}\mid {\cF^{r}}\bigg]\nonumber\\
	&+(1-p_i)\bigg(1+\frac{1}{\epsilon_i}\bigg)\|\nabla f_{\mu,i}(x^{r})-\nabla f_{\mu,i}(y_i^{r-1})\|^2,\label{eq:del:fmui}
	\end{align} where in ${\rm (i)}$ we first apply \eqref{eq:rel3}.  Note that  when $\cF^r$ is given the randomness of the first and second term in $\eqref{eq:del:fmui0}$ come from $x^{r+1}$ and $y_i^r$ respectively. Therefore, equality ${\rm (ii)}$ is true because $y_i^r=	x^r,\text{~with probability~} p_{i}$, and $y_i^r=y_{i}^{r-1}, \text{~with probability~} 1-p_{i}.$	
Setting $\epsilon_i=\frac{2}{p_i}$, the second part of \eqref{eq:p:diff} can be bounded as
\begin{align}
&\mathbb{E}_{\cJ^r}\bigg[\sum_{i=1}^N\frac{4}{p_i\rho_i}\bigg\|\nabla f_{\mu,i}(x^{r+1})-\nabla f_{\mu,i}(y_i^{r})\bigg\|^2 -
 \frac{4}{p_i\rho_i}\bigg\|\nabla f_{\mu,i}(x^{r})-\nabla f_{\mu,i}(y_i^{r-1})\bigg\|^2\mid {\cF^{r}}\bigg]\nonumber\\
&\le \sum_{i=1}^N\frac{4L_{\mu,i}^2(2+p_i)}{p_i^2\rho_i}\mathbb{E}_{\cJ^r}\|x^{r+1}-x^{r}\|^2\nonumber\\
&-\sum_{i=1}^N\frac{4(1+p_i)}{2\rho_i}\bigg\| \nabla f_{\mu,i}(x^{r})-\nabla f_{\mu,i}(y_i^{r-1})\bigg\|^2\label{eq:diff:p:2}.
\end{align}
{\bf Step 3).} In this step we combine the results from the previous steps to obtain the desired descent estimate. Combining \eqref{eq:diff:p:1} and \eqref{eq:diff:p:2} eventually we have 
\begin{align}
&\mathbb{E}_{\cJ^r}[\tilde{Q}^{r+1}- \tilde{Q}^{r}\mid \cF^r]\nonumber\\
&\le  \sum_{i=1}^N\bigg(\frac{3\beta}{2p_i} -\frac{4(1+p_i)}{2\rho_i}\bigg)\bigg\|\nabla f_{\mu,i}({  x^{r}})-\nabla f_{\mu,i}({  y_i^{r-1}})\bigg\|^2\nonumber\\
& +\sum_{i=1}^N\bigg(\frac{4L_{\mu,i}^2(2+p_i)}{p_i^2\rho_i}+\frac{L_{\mu,i}}{2}-\frac{\rho_i}{2}\bigg)\mathbb{E}_{\cJ^r}\bigg[\|x^{r+1}-x^{r}\|^2\mid\cF^r\bigg]\nonumber\\  
&+ \sum_{i=1}^{N}\frac{3\beta}{2p_i}\bigg\|\nabla f_{\mu,i}(y_i^{r-1})- \bar{G}_{\mu,i}(y^{r-1}_{i}, \phi^{(r-1)(i)},\xi^{(r-1)(i)})\bigg\|^2 
+\frac{3\tilde{p}\beta\tilde{\sigma}^2}{2J}. \label{eq:cond:Q}
\end{align}
Using the properties of conditional expectation we have
\begin{align}
&\mathbb{E}[\tilde{Q}^{r+1}- \tilde{Q}^{r}] = \mathbb{E}_{\cF^r}\bigg[\mathbb{E}_{\cJ^r}[\tilde{Q}^{r+1}- \tilde{Q}^{r}\mid \cF^r]\bigg].
\end{align}
Plugging \eqref{eq:cond:Q} in this relationship and utilizing \eqref{eq:Q:r-1}, and the definition of $\beta:=1/\sum_{i=1}^{N}\rho_{i}$, yield
\begin{align}
&\mathbb{E}[\tilde{Q}^{r+1}- \tilde{Q}^{r}]\nonumber\\
&\le
\sum_{i=1}^N\bigg(\frac{3\beta}{2p_i} -\frac{4(1+p_i)}{2\rho_i}\bigg)\mathbb{E}\bigg\|\nabla f_{\mu,i}({  x^{r}})-\nabla f_{\mu,i}({  y_i^{r-1}})\bigg\|^2\nonumber\\
&+\sum_{i=1}^N\bigg(\frac{4L_{\mu,i}^2(2+p_i)}{p_i^2\rho_i}+\frac{L_{\mu,i}}{2}-\frac{\rho_i}{2}\bigg)\mathbb{E}\bigg[\|x^{r+1}-x^{r}\|^2\bigg]
+\frac{3\tilde{p}\beta\tilde{\sigma}^2}{J}.
\end{align}
Let us define $\{\tilde{c}_i\}$ and $\hat{c}$ as following
\begin{align*}
\tilde{c}_i &= \frac{3\beta}{2p_i} -\frac{4(1+p_i)}{2\rho_i},\quad \hat{c} = \sum_{i=1}^N\bigg(\frac{4L_{\mu,i}^2(2+p_i)}{p_i^2\rho_i}+\frac{L_{\mu,i}}{2}-\frac{\rho_i}{2}\bigg).
\end{align*}
In order to prove the lemma it remains to prove that  $\tilde{c}_i<-\frac{1}{2\rho_i} \; \forall~i$, and $\hat{c}<-\sum_{i=1}^N\frac{\rho_i}{100}$. If we set $p_i = \frac{\rho_i}{\sum_{i=1}^N\rho_i}$, then we have the following
\begin{align*}
\tilde{c}_i = \frac{3}{2\rho_i}-\frac{4(1+p_i)}{2\rho_i}\leq \frac{3}{2\rho_i}-\frac{4}{2\rho_i}=-\frac{1}{2\rho_i}.
\end{align*}
To show that $\hat{c}\le -\sum_{i=1}^N\frac{\rho_i}{100}$, it is sufficient to have 
\begin{align}
\frac{4L_{\mu,i}^2(2+p_i)}{p_i^2\rho_i}+\frac{L_{\mu,i}}{2}-\frac{\rho_i}{2}\le -\frac{\rho_i}{100}.
\end{align}
It is easy to check that this inequality holds true for $\rho_i\ge\frac{5.5L_{\mu,i}}{p_i}$. 
The lemma is proved. \QED

\noindent\subsection{\bf Proof of Theorem \ref{thm:rate:zonestt-g}}
	Here we only prove the first part of the theorem. Similar steps can be followed to prove the second part. First let us define the smoothed version of optimality gap  as follows
		\begin{align}
		\Psi_\mu^r = \frac{1}{\beta^2}\mathbb{E}\bigg\|x^r-\prox_h^{1/\beta}[x^r-\beta \nabla f_\mu(x^r)]\bigg\|^2.
		\end{align}
	 We bound the gap in the following way
		\begin{align}\label{eq:Q:bound}
		&\frac{1}{\beta^2}\big\|x^r-\prox_h^{1/\beta}[x^r-\beta \nabla f_{\mu}(x^r)]\big\|^2\nonumber\\
		&\stackrel{\rm(i)}= \frac{1}{\beta^2}\big\|x^{r}-x^{r+1}+\prox_h^{1/\beta}(u^{r+1})-\prox_h^{1/\beta}[x^r-\beta \nabla f_{\mu}(x^r) ]\big\|^2\nonumber\\
		&\stackrel{\rm(ii)}\leq \frac{2}{\beta^2}\|x^{r+1}-x^r\|^2+\frac{2}{\beta^2}\|\beta \nabla f_{\mu}(x^r)+u^{r+1}-x^r\|^2\nonumber\\
		&\stackrel{\eqref{eq:u}}= \frac{2}{\beta^2}\|x^{r+1}-x^r\|^2 + 2\|\sum_{i=1}^{N} \nabla f_{\mu,i}(x^r)-v_{i_r}^r\|^2,
		\end{align} where $\rm(i)$ is true  due to \eqref{eq:x:prox}; $\rm(ii)$ is true due to the nonexpansivness of the prox operator, and equation \eqref{eq:rel2}. 
	Taking expectation on both sides yields
		\begin{align}\label{eq:Q:bound2}
		&\Psi_\mu^r\le \frac{2}{\beta^2}\mathbb{E}\|x^{r+1}-x^r\|^2 +2\mathbb{E}\|\sum_{i=1}^{N} \nabla f_{\mu,i}(x^r)-v_{i_r}^r\|^2\nonumber\\
		&\stackrel{\rm(i)}\leq \frac{2}{\beta^2}\mathbb{E}\|x^{r+1}-x^r\|^2+\frac{6}{\beta}\sum_{i=1}^N\frac{1}{\rho_i}\mathbb{E}\bigg\|\nabla f_{\mu,i}(x^r)-\nabla f_{\mu,i}(y_i^{r-1})\bigg\|^2+\frac{12\tilde{p}\tilde{\sigma}^2}{J}\nonumber\\
		&\stackrel{\eqref{eq:zone-ns:descent}}\le \frac{200}{\beta}\mathbb{E}[\tilde{Q}^r-\tilde{Q}^{r+1} ]+\frac{612\tilde{p}\tilde{\sigma}^2}{J}\nonumber\\ 
		&\stackrel{\rm(ii)}= {  1100}\bigg(\sum_{i=1}^N\sqrt{L_{\mu,i}}\bigg)^2\mathbb{E}[\tilde{Q}^r-\tilde{Q}^{r+1} ]+\frac{612\tilde{p}\tilde{\sigma}^2}{J},
		\end{align}  where in $\rm(i)$ we utilize \eqref{eq:vir_bd}. To get $\rm(ii)$ let us pick $\rho_{i}=\frac{5.5L_{\mu,i}}{p_i}$, therefore we have
	$\rho_{i}=5.5L_{\mu,i}\frac{\sum_{i=1}^{N}\rho_i}{\rho_{i}},$
	which leads to $\rho_i=\sqrt{5.5L_{\mu,i}\sum_{j=1}^{N}\rho_j}=\sqrt{5.5L_{\mu,i}}\sqrt{\sum_{j=1}^{N}\rho_j}.$ Summing both sides over $i=1,2,\cdots N$, and simplifying the result we get 
	$$\sqrt{\sum_{i=1}^{N}\rho_i}=\sum_{i=1}^{N}\sqrt{5.5L_{\mu,i}}.$$ Finally, squaring both sides and  set $\beta:=1/\sum_{i=1}^{N}\rho_i$ we reach $\frac{1}{\beta}=5.5(\sum_{i=1}^{N}\sqrt{L_{\mu,i}})^2.$
	Let us sum both sides of \eqref{eq:Q:bound2} over $T$ iterations, use telescopic property, and divide both sides by $T$, we obtain
		\begin{align*}
		\frac{1}{T}\sum_{r=1}^T\Psi_\mu^r\leq 1100\bigg(\sum_{i=1}^N\sqrt{L_{\mu,i}}\bigg)^2\frac{\mathbb{E}[\tilde{Q}^1-\tilde{Q}^{T+1}]}{T} + \frac{612\tilde{p}\tilde{\sigma}^2}{J}.
		\end{align*}
	Since $u$ is uniformly random number in $\{1,2,\cdots, T\}$, we finally have
			\begin{align}
			\mathbb{E}_{u}[\Psi_\mu^u]\leq {  1100}\Big(\sum_{i=1}^N\sqrt{L_{\mu,i}}\Big)^2\frac{\mathbb{E}[{\tilde{Q}^1-\tilde{Q}^{T+1}} ]}{T}+ \frac{612\tilde{p}\tilde{\sigma}^2}{J}.\label{eq:sublin}
			\end{align}
	
	Now let us bound the gap $\Psi^r$. Using the definition of $\Psi^r$ we have
		\begin{align*}
		\Psi^r&=\frac{1}{\beta^2}\mathbb{E}\bigg[\|x^r-\prox_h^{1/\beta}[x^r-\beta \nabla f(x^r)]\bigg]\nonumber\\
		&=\frac{1}{\beta^2}\mathbb{E}\bigg[\|x^r-\prox_h^{1/\beta}[x^r-\beta \nabla f(x^r)]\nonumber\\
		&-\prox_h^{1/\beta}[x^r-\beta \nabla f_\mu(x^r)]+\prox_h^{1/\beta}[x^r-\beta \nabla f_\mu(x^r)]\|^2\bigg]\\
		&\stackrel{\rm(i)}\le 2\Psi_\mu^r + \frac{\mu^2{L}^2(M+3)^3}{2},
		\end{align*} where in $\rm(i)$ we use  $\eqref{eq:rel2}$; the nonexpansiveness of the prox operator; and inequality \eqref{eq:diff:grad:bd}. Next because $r$ is a uniformly random number picked form $\{1,2,\cdots, T\}$ we have
		\begin{align}
			&\mathbb{E}_{u}[\Psi^u]\leq 2 \mathbb{E}_{u}[\Psi_\mu^u]  + \frac{\mu^2{L}^2(M+3)^3}{2}\nonumber\\
			&\stackrel{\eqref{eq:sublin}}\le {  2200}\Big(\sum_{i=1}^N\sqrt{L_{\mu,i}}\Big)^2\frac{\mathbb{E}[{\tilde{Q}^1-\tilde{Q}^{T+1}}]}{T} +\frac{\mu^2{L}^2(M+3)^3}{2} + \frac{1024\tilde{p}\tilde{\sigma}^2}{J}.
		\end{align}
		The proof is complete. \QED
{
	\bibliography{ref_davood}
}

\end{document}